\documentclass{amsart}
\usepackage{amsmath}%
\usepackage{amsfonts}%
\usepackage{amssymb}%
\usepackage{graphicx,pb-diagram,color}
\newtheorem{theorem}{Theorem}
\theoremstyle{plain}

\newtheorem{corollary}[theorem]{Corollary}

\newtheorem{definition}[theorem]{Definition}
\newtheorem{example}[theorem]{Example}

\newtheorem{problem}{Problem}

\newtheorem{remark}[theorem]{Remark}

\numberwithin{equation}{section}
\def\C{{\mathbb C}}

\def\M{{\mathrm M}}
\def\EF{\mathcal{E}}

\def\R{{\mathbb R}}
\def\Z{{\mathbb Z}}

\def \P{{\mathbb P}}

\def\S{\mathbb{S}}

\def\CD{{\mathcal D}}
\def\CH{{\mathcal H}}
\def\CM{{\mathcal M}}
\def\CU{{\mathcal U}}
\def\CV{{\mathcal V}}
\def\CW{{\mathcal W}}
\def\tmu{\widetilde{\mu}}

\def\CHd{\mathcal{H}_{(d)}}
\def\PHd{\P(\mathcal{H}_{(d)})}

\def\ra{\rightarrow}
\def\Pn{{\P(\C^{n+1})}}

\def\XXint#1#2#3{{\setbox0=\hbox{$#1{#2#3}{\int}$}
\vcenter{\hbox{$#2#3$}}\kern-.5\wd0}}

\begin{document}

\title[Polynomial System Solving]{The complexity and geometry of numerically solving polynomial systems.}


\author{Carlos Beltr\'an}
\thanks{C. B.: depto. de Matem\'aticas, Estad\'istica y Computaci\'on, Universidad de Cantabria, Santander, Spain. Research was Partially Supported by MTM2010-16051, Spanish Ministry of Science. {\tt carlos.beltran@unican.es}}

\author{Michael Shub}
\thanks{M.S.: CONICET, IMAS, Universidad de Buenos Aires, Argentina and CUNY Graduate School, New York, NY, USA.
M.S. was partially
supported by a CONICET grant PIP0801 2010--2012 and by ANPCyT PICT 2010--00681. {\tt shub.michael@gmail.com}}


\date{\today}

\dedicatory{This paper is dedicated to our beloved friend and colleague Jean Pierre Dedieu.}

\begin{abstract}
These pages contain a short overview on the state of the art of efficient numerical analysis methods that solve systems of
multivariate polynomial equations. We focus on the work of Steve
Smale who initiated this research framework, and on the collaboration
between Stephen Smale and Michael Shub, which set the foundations of this
approach to polynomial system--solving, culminating in the more recent advances of Carlos Beltr\'an, Luis Miguel Pardo, Peter B\"urgisser and Felipe Cucker.
\end{abstract}

\maketitle
\section{The modern numerical approach to polynomial system solving}\label{section:1}
In this paper we survey some of the recent advances in the solution of polynomial systems. Such a classical topic has been studied by hundreds of authors from many different perspectives. We do not intend to make a complete historical description of all the advances achieved during the last century or two, but rather to describe in some detail the state of the art of what we think is the most successful (both from practical and theoretical perspectives) approach. { Homotopy methods} are used to solve polynomial systems in real life applications all around the world.

The key ingredient of homotopy methods is a one--line thought: given a goal system to be solved, choose some other system (similar in form, say with the same degree and number of variables) with a known solution $\zeta_0$, and move this new system to the goal system, tracking how the known solution moves to a solution of the goal. Before stating any notation, we can explain briefly why this process is reasonable: if for every $t\in[0,1]$ we have a system of equations $f_t$ ($f_0$ is the system with a known solution, $f_1$ is the one we want to solve), then we are looking for a path $\zeta_t$, $t\in[0,1]$, such that $f_t(\zeta_t)=0$. As long as the derivative $df_t(\zeta_t)$ is invertible for all $t$ we can continue the solution from $f_0$ to $f_1$, by the implicit function theorem. Now we have various methods to accomplish this continuation. We can slowly increment $t$ and use iterative numerical solution methods such as Newton's method to track the solution or we may differentiate the expresion $f_t(\zeta_t)=0$ and solve for $d/(dt)(\zeta_t)=\dot\zeta_t$. Then, we can write our problem as an initial value problem:
\begin{equation}\label{eq:EDO}
\begin{cases}\dot\zeta_t=Df_t(\zeta_t)^{-1}f_t(\zeta_t)\\ \zeta_0\text{ known}\end{cases}
\end{equation}
Systems of ODEs have been much studied and hence this is an interesting idea: we have reduced our original problem to a very much studied one. One can just plug in a standard numerical ODE solver such as backward Euler or a version of Runge--Kutta's method. Even then, in practice, it is desirable to, from time to time, perform some steps of Newton's method $z\rightarrow x-Df_t(x)^{-1}f_t(x)$ to our approximation $z_t$ of $\zeta_t$, to get closer to the path $(f_t,\zeta_t)$. After some testing and adjustment of parameters, this na\"ive idea can be made to work with impressive practical performance and there are several software packages which attain spectacular results (solving systems with many variables and high degree) in a surprisingly short running time, see for example \cite{HOM4PSwww,Bertini,V99,Leykin:NAG4M2}

From a mathematical point of view, there are several things in the process we have just described that need to be analyzed: will there actually exist a path $\zeta_t$ (maybe it is only defined for, say, $t<1/2$)? what is the expected complexity of the process (in particular, can we expect average polynomial running time in some sense)? what ``simple system with a known solution'' should we start at? how should we join $f_0$ and $f_1$, that is what should be the path $f_t$?

In the last few decades a lot of progress has been made in studying these questions. This progress is the topic of this paper.

\section{A technical description of the problem}
We will center our attention in Smale's 17--th problem, which we recall now.
\begin{problem}\label{prob:17}
Can a zero of n complex polynomial equations in n unknowns be { found approximately},  { on the average, in polynomial time} with a { uniform
algorithm}?
\end{problem}
We have written in bold the technical terms that need to be clarified.

In order to understand the details of the problem and the solution suggested in Section  \ref{section:1}, we need to describe some important concepts and notation in detail. Maybe the first one is our understanding of what a ``solution'' is: clearly, one cannot expect solutions of polynomial systems to be rational numbers, so one can only search for ``quasi--solutions'' in some sense. There are several definitions of such a thing, the most stable being the following one (introduced in \cite{Smale81}, see also \cite{Kim85,Kim88}):
\begin{definition}\label{def:appzero}
Given a polynomial system, understood as a mapping $f:\C^n\rightarrow\C^n$, an approximate zero of $f$ with associated (exact) zero $\zeta$ is a vector $z_0\in\C^n$ such that
\[
\|z_k-\zeta\|\leq\frac{1}{2^{2^k-1}}\|z_0-\zeta\|,\quad k\geq0,
\]
where $z_k$ is the result of applying $k$ times Newton's operator $z\mapsto z-Df(z)^{-1}f(z)$ (note that the definition of approximate zero implicitly assumes that $z_k$ is defined for all $k\geq0$.)
\end{definition}
The power of this definition is that, as we will see below, given any polynomial system $f$ and any exact zero $\zeta\in\C^n$, approximate zeros of $f$ with associated zero $\zeta$ exist whenever $Df(\zeta)$ is an invertible matrix.

Recall that our first goal is to transform the problem of polynomial system solving into an implicit function problem or an ODE system like that of (\ref{eq:EDO}). There exist two principal reasons why the solution of such a system can fail to be defined for all $t>0$: that the function defining the derivative is not everywhere defined (this corresponds naturally to $Df_t(\zeta_t)$ not being invertible), and that the solution escapes to infinity. The first problem seems to be more delicate and difficult to solve, but the second one is actually very easily dealt with: we just need to define our ODE in a compact manifold, instead of just in $\C^n$. The most similar compact manifold to $\C^n$ is $\Pn$, and the way to take the problem into $\Pn$ is just homogenizing the equations.

\begin{definition}
Let $f:\C^n \rightarrow \C^n$ be a polynomial system, that is $f=(f_1,\ldots,f_n)$ where $f_i:\C^n\rightarrow\C$ is a polynomial of degree some $d_i$,
\[
f(x_1,\ldots,x_n)=\sum_{\alpha_1+\cdots+\alpha_n\leq d_i}a^{(i)}_{\alpha_1,\ldots,\alpha_n}x_1^{\alpha_1}\cdots x_n^{\alpha_n}.
\]
The homogeneous counterpart of $f$ is $h:\C^{n+1}\rightarrow\C^n$ defined by $h=(h_1,\ldots,h_n)$ where
\[
h(x_0,x_1,\ldots,x_n)=\sum_{\alpha_1+\cdots+\alpha_n\leq d_i}a^{(i)}_{\alpha_1,\ldots,\alpha_n}x_0^{d_i-\sum_{i=1}^n\alpha_i}x_1^{\alpha_1}\cdots x_n^{\alpha_n}.
\]
We will talk about such a system $h$ simply as an { homogeneous system}.
\end{definition}
Note that if $\zeta$ is a zero of $f$ then $(1,\zeta)$ is a zero of the homogeneous counterpart $h$ of $f$. Reciprocally, if $\zeta=(\zeta_0,\zeta_1,\ldots,\zeta_n)$ is a zero of $h$ and if $\zeta_0\neq0$, then $(\zeta_1/\zeta_0,\ldots,\zeta_n/\zeta_0)$ is a zero of $f$. Thus, the zeros of $f$ and $h$ are in correspondence and we can think of solving $h$ and then recovering the zeros of $f$ (this is not a completely obvious process when we only have approximate zeros, see \cite{BePa:JAMS}.) Moreover, it is clear that for any complex number $\lambda\in\C$ and for $x\in\C^{n+1}$ we have
\[
h(\lambda x)=Diag(\lambda^{d_1},\ldots,\lambda^{d_n})h(x),
\]
and thus the zeros of $h$ lie naturally in the projective space $\Pn$.

As we will be working with homogeneous systems and projective zeros, we need a definition of approximate zero in the spirit of Definition \ref{def:appzero} which is amenable to a projective setting. The following one, which uses the projective version \cite{Sh93} of Newton's operator, makes the work. Here and throughout the paper, given a matrix or vector $A$, by $A^*$ we mean the complex conjugate transpose of $A$, and by $d_R(x,y)$ we mean the Riemannian distance from $x$ to $y$, where $x$ and $y$ are elements in some Riemannian manifold.
\begin{definition}\label{def:appzeroproj}
Given an homogeneous system $h$, an approximate zero of $h$ with associated (exact) zero $\zeta\in\Pn$ is a vector $z_0\in\Pn$ such that
\[
d_R(z_k,\zeta)\leq\frac{1}{2^{2^k-1}}d_R(z_0,\zeta),\quad k\geq0,
\]
where $z_k$ is the result of applying $k$ times the projective Newton operator $z\mapsto z-Dh(z)\mid_{z^\perp}^{-1}h(z)$ (again, the definition of approximate zero implicitly assumes that $z_k$ is defined for all $k\geq0$.) Here, by $Df(z)\mid_{z^\perp}$ we mean the restriction of the derivative of $h$ at $z$, to the (complex) orthogonal subspace $z^\perp=\{y\in\C^{n+1}:y^*x=0\}$.
\end{definition}
It is a simple exercise to verify that (projective) Newton's method is well defined, that is the point it defines in projective space does not depend on the representative $z\in\C^{n+1}$ chosen for a point in projective space.

A (projective) approximate zero of $h$ is thus a projective point such that the successive iterates of the projective Newton operator quickly approach an exact zero of $h$. Thus finding an approximate zero is an excellent output of a numerical zero--finding algorithm to solve $h$.

Because we are going to consider paths of systems $\{h_t\}_{t\in[a,b]}$, it is convenient to fix a framework where one can define these nicely. To this end, we consider the vector space of homogeneous polynomials of fixed degree $s\geq1$:
\[
 \CH_s=\{h\in\C[x_0,\ldots,x_n]:\,h\text{ is homogeneous of degree }s\}.
\]
It is convenient to consider an Hermitian product (and the associated metric) on $\CH_s$. A desirable property of such a metric is the unitary invariance, namely, we would like to have an Hermitian product such that
\[
\langle h,g\rangle_{\CH_s}=\langle h\circ U,g\circ U\rangle_{\CH_s},\quad \forall\;U\in\CU_{n+1},
\]
where $\CU_{n+1}$ is the group of unitary matrices of size $n+1$. Such property was studied in detail in \cite{ShSm:BezI}. It turns out that there exists a unique (up to scalar multiplication) Hermitian product that satisfies it, the one defined as follows:
\[
\langle \sum_{\alpha_0+\cdots+\alpha_n= s}a^{(i)}_{\alpha_1,\ldots,\alpha_n}x_0^{\alpha_0}\cdots x_n^{\alpha_n},\sum_{\alpha_0+\cdots+\alpha_n=s}b^{(i)}_{\alpha_1,\ldots,\alpha_n}x_0^{\alpha_0}\cdots x_n^{\alpha_n}\rangle_{\CH_s}=
\]
\[
\sum_{\alpha_0+\cdots+\alpha_n= s}\frac{\alpha_0!\cdots\alpha_n!}{s!}a^{(i)}_{\alpha_1,\ldots,\alpha_n}\overline{b^{(i)}_{\alpha_1,\ldots,\alpha_n}},
\]
where $\overline{\cdot}$ just means complex conjugation. Note that this is just a weighted version of the standard complex Hermitian product in complex affine space.

Then, given a list of degrees $(d)=(d_1,\ldots,d_n)$, we consider the vector space
\[
 \CHd=\prod_{i=1}^n\CH_{d_i}.
\]
Note that an element $h$ of $\CHd$ can be seen both as a mapping $h:\C^{n+1}\rightarrow\C^n$ or as a polynomial system, and can be identified by the list of coefficients of $h_1,\ldots,h_n$. We denote by $\PHd$ the projective space associated to $\CHd$, by $N$ the complex dimension of $\PHd$ (so the dimension of $\CHd$ is $N+1$) and we consider the following Hermitian structure in $\CHd$:
\[
\langle h,g\rangle=\sum_{i=1}^n\langle h_i,g_i\rangle_{H_{d_i}},\quad \|h\|=\langle h,h\rangle^{1/2}.
\]
This Hermitian product (and the associate Hermitian structure and metric) is also called the Bombieri--Weyl or the Kostlan product (structure, metric). As usual, this Hermitian product in $\CHd$ defines an associated Riemannian structure given by the real part of $\langle\cdot,\cdot\rangle$. We can thus consider integrals of functions defined on $\CHd$.

We denote by $\S$ the unit sphere in $\CHd$, and we endow $\S$ with the inherited Riemannian structure from that of $\CHd$. Then, $\PHd$ has a natural Riemannian structure, the unique one making the projection $\S\rightarrow\PHd$ a Riemannian submersion. That is the derivative of the projection restricted to the normal to the fibers is an isometry. We can thus also consider integrals of functions defined in $\S$ or $\PHd$. We can now talk about probabilities in $\S$ or $\PHd$: given a measurable (nonnegative or integrable) mapping $X$ defined in $\S$ or $\PHd$, we can consider its expected value:
\[
{\rm E}_\S(X)=\frac{1}{\nu(\S)}\int_\S X(h)\,dh \quad \text{or}\quad {\rm E}_{\PHd}(X)=\frac{1}{\nu(\PHd)}\int_{\PHd} X(h)\,dh,
\]
where we simply denote by $\nu(E)$ the volume of a Riemannian manifold $E$. Similarly, one can talk about probabilities in $\CHd$ according to the standard Gaussian distribution compatible with $\langle\cdot,\cdot\rangle$: given a measurable (nonnegative or integrable) mapping $X$ defined in $\CHd$, its expected value is:
\[
{\rm E}_{\CHd}(X)=\frac{1}{(2\pi)^{N+1}}\int_{\CHd}X(h)e^{-\|h\|^2/2}\,dh.
\]
We can now come back to Problem \ref{prob:17} and see what do each of the terms in that problem mean: Smale himself points out that one can just solve homogeneous systems (as suggested above). We still have a few terms to clarify:
\begin{itemize}
\item { found approximately.} This means finding an approximate zero in the sense of Definition \ref{def:appzeroproj}.
\item { on the average, in polynomial time.} This now means that, if $X(h)$ is the time needed by the algorithm to output an approximate zero of the input system $h$, then the expected value of $X$ is a quantity polynomial in the input size, that is polynomial in $N$. The number of variables, $n$, and the maximum of the degrees, $d$, are smaller than $N$, and hence one attempts to get a bound on the expected value of $X$, as a polynomial in $n,d,N$.
\item { uniform algorithm.} Smale demands an algorithm in the Blum--Shub--Smale model \cite{BlShSm,BlCuShSm98}, that is exact operations and comparisons between real numbers are assumed. This assumption departs from the actual performance of our computers, but it is close enough to be translated to performance in many situations. Uniform means that the same algorithm works for all $(d)$ and $n$.
\end{itemize}

\section{Geometry and condition number}\label{sec:geometry}
We can now set up a geometric framework for homotopy methods. Consider the following set, usually called the { solution variety}:
\begin{equation}\label{eq:V}
\CV=\{(h,\zeta)\in\PHd\times\Pn:h(\zeta)=0\}.
\end{equation}
This set is actually a smooth complex submanifold (as well as a complex algebraic subvariety) of $\PHd\times\Pn$, see \cite{BlCuShSm98}, and is clearly compact. It will be useful to consider the following diagram.
\begin{equation}\label{eq:pi1pi2}
\begin{matrix}
&&\CV&&\\
&\pi_1\swarrow&&\searrow\pi_2&\\
&&&&\\
\PHd&&&&\Pn
\end{matrix}
\end{equation}
It is clear that $\pi_1^{-1}(h)$ is a copy of the zero set of $h$. Reciprocally, for fixed $\zeta\in\Pn$, the set $\pi_2^{-1}(\zeta)$ is the vector space of polynomial systems that have $\zeta$ as a zero.

Let $\Sigma'\subseteq\CV$ be the set of critical points of $\pi_1$ and $\Sigma=\pi_1(\Sigma')\subseteq\PHd$ the set of critical values of $\pi_1$. It is not hard to prove that:
\begin{itemize}
\item  $\pi_1$ restricted to the set $\CV\setminus\pi_1^{-1}(\Sigma)$ is a (smooth) $\CD$--fold covering map, where $\CD=d_1\cdots d_n$ is the Bez\'out number.
\item $\Sigma'=\{(h,\zeta)\in\CV:Dh(\zeta)\mid_{\zeta^{\perp}}\text{ has non--maximal rank}\}$. In that case, we say that $\zeta$ is a { singular zero} of $h$. Otherwise, we say that $\zeta$ is a { regular zero} of $h$.
\end{itemize}
This means, in particular, that the homotopy process described above can be carried out whenever the path of systems lies outside of $\Sigma$:
\begin{theorem}\label{th:lift}
Let $\{h_t:t\in[a,b]\}$ be a $C^1$ curve in $\PHd\setminus\Sigma$ and let $\zeta$ be a zero of $h_a$. Then, there exists a unique lift of $h_t$ through $\pi_1$, that is a $C^1$ curve $(h_t,\zeta_t)\in\CV$ such that $\zeta_a=\zeta$. In particular, $\zeta_b$ is a zero of $h_b$. Moreover, the lifted curve satisfies:
\begin{equation}\label{eq:ecdif}
\frac{d}{dt}(h_t,\zeta_t)=\left(\dot h_t,Dh_t(\zeta_t)\mid_{\zeta_t^{\perp}}^{-1}\dot h_t(\zeta_t)\right).
\end{equation}
Finally, the set $\Sigma\subseteq\PHd$ is a complex proper algebraic variety, thus it has real codimension $2$ and the projection of most real lines in $\CHd$ to $\PHd$ does not intersect $\Sigma$.
\end{theorem}
In the case the thesis of Theorem \ref{th:lift} holds we just say that $\zeta_a$ can be continued to a zero $\zeta_b$ of $h_a$. One can be even more precise:
\begin{theorem}\label{th:lift2}
Let $\{h_t:t\in[a,b]\}$ be a $C^1$ curve in $\PHd\setminus\Sigma$ and let $\zeta$ be a zero of $h_a$. Then, every zero $\zeta$ of $h_a$ can be continued to a zero of $h_b$, defining a bijection between the $\CD$ zeros of $h_a$ and those of $h_b$.
\end{theorem}
\begin{remark}
Even if $h_t$ cross $\Sigma$ some solutions may be able to be continued while others may not.
\end{remark}

The (normalized) condition number \cite{ShSm:BezI} is a quantity describing ``how close to singular'' a zero is. Given $h\in\CHd$ and $z\in\Pn$, let
\begin{equation}\label{eq:mun}
\mu( f,z)=\| f\|\|(Dh(z)\mid_{z^\perp})^{-1}Diag(\|z\|^{d_i-1}d_i^{1/2})\|_2,
\end{equation}
and $\mu(f,z)=+\infty$ if $Dh(z)\mid_{z^\perp}$ is not invertible. Sometimes $\mu$ is denoted $\mu_{\rm norm}$ or $\mu_{\rm proj}$ but we prefer to keep the more simple notation here. One of the most important properties of $\mu$ is that it is an upper bound for the norm of the (locally defined) implicit function related to $\pi_1$ in (\ref{eq:pi1pi2}). Namely, let $(\dot h,\dot \zeta)\in T_{(h,\zeta)}\CV$ where $(h,\zeta)\in\CV$ is such that $\mu(h,\zeta)<+\infty$. Then,
\begin{equation}\label{eq:muboundszetadot}
\|\dot\zeta\|\leq \mu(h,\zeta)\|\dot h\|,\quad \mu(h,\zeta)\geq\sqrt{n}.
\end{equation}

We also have the following result.
\begin{theorem}[Condition Number Theorem,\cite{ShSm:BezI}]\label{th:condnumber}
\[
\mu(h,\zeta)=\frac{1}{\sin\left(d_R(h,\Sigma_{\zeta})\right)},
\]
where $d_R$ is the Riemannian distance in $\P(\CHd)$ and
\[
\Sigma_{\zeta}=\{h\in\P(\CHd)\; :\; h(\zeta)=0,\; {\hbox {\rm
and}}\; Dh(\zeta)\mid_{\zeta^\perp} {\hbox {\rm\; is not invertible}}\}.
\]
\end{theorem}
Note that this is a version of the classical Condition Number Theorem of linear algebra (see Theorem \ref{th:cnt} below). The existence of approximate zeros in the sense of Definition \ref{def:appzeroproj} above is also guaranteed by this condition number, as was noted in \cite{ShSm:BezI}. More precisely:

\begin{theorem}[$\mu$--Theorem, \cite{ShSm:BezI}]\label{mu:teor}
There exists a constant $u_0>0$ ($u_0=0.17586$ suffices) with the following property. Let $(h,\zeta)\in \CV$ and let $z\in\Pn$ satisfy
\[
d_R(z,\zeta)\leq\frac{u_0}{d^{3/2}\,\mu(h,\zeta)}.
\]
Then, $z$ is an approximate zero of $h$ with associated zero $\zeta$.
\end{theorem}

\section{The complexity of following a homotopy path}\label{sec:c1hp}
The sentence ``can be continued'' in the discussion of Section \ref{sec:geometry} can be made much more precise, by defining an actual path--following method. It turns out that the unique method that has actually been proved to correctly follow the homotopy paths and at the same time achieve some known complexity bound is the most simple one, which only uses the projective Newton operator, and not an ODE solver step.
\begin{problem}
It would be an interesting project to compare the overall cost of using a higher order ODE solver to the projective Newton--based method we describe below. Higher order methods or even predictor--corrector methods may require fewer steps but be more expensive at each step so a total cost comparison is in order. Some experience indicates that higher order methods are rarely cheaper, if ever. See \cite{BorodinMunro1975,Kim85,Kim88}.
\end{problem}
More precisely, the projective Newton--based homotopy method is as follows. Given a $C^1$ path $\{h_t:a\leq t\leq b\}\subseteq\PHd$, and given $z_a$ an approximate zero of $h_a$ with associated (exact) zero $\zeta_a$, let $t_0>0$ be ``small enough'' and let
\[
z_{a+t_0}=z_a-(Dh_{a+t_0}(z_a)\mid_{z_a^{\perp}})^{-1}h_{a+t_0}(z_a),
\]
that is $z_{a+t_0}$ is the result of one application of the projective Newton operator based on $h_{a+t_0}$ to the point $z_a$. If $z_a$ is an approximate zero of $h_a$ and $t_0$ is small enough, then $z_a$ can be close enough to the actual zero $\zeta_{a+t_0}$ of $h_{a+t_0}$ to satisfy Theorem \ref{mu:teor} and thus be an approximate zero of $h_{a+t_0}$ as well. Then, by definition of approximate zero, $z_{a+t_0}$ will be half--closer to $\zeta_{a+t_0}$ than $z_a$. This leads to an inductive process (choosing $t_1$, then $t_2$, etc. until $h_b$ is reached) that, analysed in detail, can be made to work and actually programmed. The details on how to choose $t_0$ would take us too far apart from the topic, so we just give an intuitive explanation: if we are to move from $(h_a,\zeta_a)$ to $(h_{a+t_0},\zeta_{a+t_0})$ we must be sure that we are far enough from $\Sigma'$ to have our algorithm behaving properly. As the condition number essentially measures the distance to $\Sigma'$, it should be clear that the bigger the condition number, the smaller step $t_0$ we can take. This idea lead to the following result (see \cite{ShSm:BezV} for a weaker, earlier result):
\begin{theorem}[\cite{Shub2007}]\label{th:BezVI}
Let $(h_t,\zeta_t)\subseteq\CV\setminus\Sigma'$, $t\in[a,b]$ be a $C^1$ path. If the steps $t_0,t_1,\ldots$ are correctly chosen, then an approximate zero of $h_b$ is reached at some point, namely there is a $k\geq1$ such that $\sum_{i=0}^kt_i=b-a$ ($k$ is the number of steps in the inductive process above.) Moreover, one can bound
\[
k\leq \lceil Cd^{3/2}L_\kappa \rceil,
\]
where $d$ is the maximum of the degrees in $(d)$, $C$ is some universal constant, and
\begin{equation}\label{eq:condmet}
L_\kappa=\int_{a}^b\mu(h_t,\zeta_t)\|(\dot h_t,\dot \zeta_t)\|\,dt
\end{equation}
is the { condition lenght} of the path $(h_t,\zeta_t)$. Moreover, the amount of arithmetic operations needed in each step is polynomial in the input size $N$, and hence the total complexity of the path--following procedure is a quantity polynomial in $N$ and linear in $L_\kappa$
\end{theorem}
There exist several ways to algorithmically produce the steps $t_0,t_1,\ldots$ in this theorem (and indeed the process has been programmed in two versions \cite{BeltranLeykin2009,BeltranLeykin2011},) but the details are too technical for this report, see \cite{Beltran2009NM,BurgisserCucker2009,DedieuMalajovichShub2011}. We also point out that, if the path we are following is linear, i.e. $h_t=(1-t)h_0+th_1$, and if the input coordinates are rational numbers, then all the operations can be carried out over the rationals without a dramatic increase of the bit size of intermediate results, see \cite{BeltranLeykin2011}.

Note that since $L_\kappa$ is a length it is independent of the $C^1$ parametrization of the path. If we specify a path of polynomial systems in $\CHd$ then we project the path of polynomials and solutions into $\CV$ to calculate the length. We may project from $\CHd$ to $\S$ first and reparametrize if we wish. For example, we project the straight line segment $h_t=(1-t)g+th$ for $0\leq t\leq 1$ into $\S$ and reparametrize by arc--length. If $\|g\|=\|h\|=1$ the resulting curve is
\[
h_t=g\cos(t)+\frac{h-\langle h,g\rangle g}{\|h-\langle h,g\rangle g\|}\sin(t)
\]
which is an arc of great circle through $g$ and $h$. If $0\leq t\leq d_R(g,h)$, then the arc goes from $g$ to $h$. Here  $d_R(g,h)$ is the Riemannian distance in $\S$ between $g$ and $h$ which is the angle between them.

\section{The problem of good starting points}\label{sec:start}
We now come back to the original question in Smale's 17-th problem. Our plan is to analyse the complexity of an algorithm that we could call ``linear homotopy'': choose some $g\in\S,\zeta\in\Pn$ such that $g(\zeta)=0$ (we will call $(g,\zeta)$ a ``starting pair''). For input $h\in\S$, consider the path contained in the great circle :
\begin{equation}\label{eq:gc}
h_t=g\cos(t)+\frac{h-\langle h,g\rangle g}{\|h-\langle h,g\rangle g\|}\sin(t),\quad t\in[0,d_R(g,h)].
\end{equation}
Then, use the method described in Theorem \ref{th:BezVI} to track how $\zeta_0$ moves to $\zeta_{d_R(g,h)}$, a zero of $h_{d_R(g,h)}=h$, thus producing an approximate zero of $h$. We call this linear homotopy (maybe a more appropriate name would be ``great circle homotopy'') because great circles are projections on $\S$ of segments in $\CHd$.

Assuming that the input $h$ is uniformly distributed on $\S$, we can give an upper bound for the average number of arithmetic operations needed for this task (that is, the average complexity of the linear homotopy method) by a polynomial in $N$ multiplied by the following quantity:
\[
\frac{1}{\nu(\S)}\int_{h\in\S}\int_0^{d_R(g,h)}\mu(h_t,\zeta_t)\|(\dot h_t,\dot \zeta_t)\|\,dt\;d\S,
\]
where $h_t$ is defined by (\ref{eq:gc}) and $\zeta_t$ is defined by continuation (the fact that $h_t\cap\Sigma=ptyset$, and thus the existence of such $\zeta_t$, is granted by Theorem \ref{th:lift} for most choices of $g,h$). It is convenient to replace this last expected value by a similar upper bound:
\[
\mathcal{A}_1(g,\zeta)=\frac{1}{\nu(\S)}\int_{h\in\S}\int_0^{\pi}\mu(h_t,\zeta_t)\|(\dot h_t,\dot \zeta_t)\|\,dt\;d\S.
\]
Note that we are just replacing the integral from $0$ to $d_R(g,h)$ by the distance from $0$ to $\pi$.

We thus have:
\begin{theorem}
Let $(g,\zeta)\in\CV$. The average complexity of linear homotopy with starting pair $(g,\zeta)$ is bounded above by a polynomial in $N$ multiplied by $\mathcal{A}_1(g,\zeta)$.
\end{theorem}
This justifies the following definition:
\begin{definition}\label{def:gsp}
Fix some polynomial\footnote{Because $n,d\leq N$, we could just talk about a one variable polynomial $p(x)$ and change $p(n,d,N)$ to $p(N)$ in the following definition. However, we prefer here to be a bit more precise.} $p\in\R[x,y,z]$. We say that $(g,\zeta)$ is a { good starting pair} w.r.t. $p(x,y,z)$ if $\mathcal{A}_1(g,\zeta)\leq p(n,d,N)$ (which implies that the average number of steps of the linear homotopy is $O(d^{3/2}p(n,d,N))$.) From now on, if nothing is said, we assume $p(x,y,z)=\sqrt{2}\pi xz$. Thus, $(g,\zeta)\in\CV$ is a good initial pair if $\mathcal{A}_1(g,\zeta)\leq \sqrt{2}\pi nN$.
\end{definition}
So, if a good initial pair is known for all choices of $n$ and the list of degrees $(d)$, then the total average complexity of linear homotopy is polynomial in $N$. In other words, finding good starting pairs for every choice of $n$ and $(d)$ gives a satisfactory solution to Problem (\ref{prob:17}).

In \cite{ShSm:BezV} the following pair\footnote{The pair conjectured in \cite{ShSm:BezV} does not contain the extra $d_i^{1/2}$ factors. There is, however, some consensus
that these extra factors should be added, for with these factors the condition number $\mu(g,\zeta)=n^{1/2}$ is
minimal.} was conjectured to be a good starting pair (for some polynomial $p(x,y)$)
:
\begin{equation}\label{eq:shsmconj}
g(z)=\begin{cases}d_1^{1/2}z_0^{d_1-1}z_1,\\\hskip .5cm \vdots
\\d_{n}^{1/2}z_0^{d_n-1}z_n\end{cases},\qquad \zeta=(1,0,\ldots,0).
\end{equation}
To this date, proving this conjecture is still an open problem. Some experimental data supporting this conjecture was shown in \cite{BeltranLeykin2009}

\subsection{Choosing initial pairs at random: an Average Las Vegas algorithm for problem (\ref{prob:17})}
One can study the average value of the quantity $\mathcal{A}_1(g,\zeta)$ described above. Most of the results in this section are based in the fact that the expected value of the square of the condition number is relatively small. This was first noted in \cite{ShSm93}, then this expected value was computed exactly in \cite{BePa:JFoCM}:
\begin{theorem}\label{th:averagemu2}
Let $h\in\S$ be chosen at random, and let $\zeta$ be chosen at random, with the uniform distribution, among the zeros of $h$. Then, the expected value of $\mu^2(h,\zeta)$ is at most $nN\CD$. More exactly:
\[
{\rm E}_{h\in\S}\left( \sum_{\zeta:h(\zeta)=0}\mu(h,\zeta)^2\right)=
\CD N\left(n\left(1+\frac{1}{n}\right)^{n+1}-2n-1\right)\leq nN\CD.
\]
In particular, in the case of one homogeneous polynomial of degree $d$ (i.e. $n=1$,) we have:
\[
{\rm E}_{h\in\S}\left( \sum_{\zeta:h(\zeta)=0}\mu(h,\zeta)^2\right)=d(d+1).
\]
\end{theorem}

Now we use some arguments which are very much inspired by ideas from integral geometry, one of the main contributions of Lluis Santal\'o to XX century mathematics. We can try to compute the expected value of $\mathcal{A}_1(g,\zeta)$. Although this can be done directly (see \cite{BeltranShub2009},) it is easier to first consider an upper bound of $\mathcal{A}_1$: let us note from (\ref{eq:muboundszetadot}) that
\begin{equation}\label{eq:A1mu2}
\mathcal{A}_1(g,\zeta)\leq  \frac{\sqrt{2}}{\nu(\S)}\int_{h\in\S}\int_0^{\pi}\mu(h_t,\zeta_t)^2\,dt\;d\S.
\end{equation}
So, we have
\[
{\rm E}_{g\in\S}\left( \sum_{\zeta:g(\zeta)=0}\mathcal{A}_1(g,\zeta)\right)\leq {\rm E}_{g\in\S}\left( \sum_{\zeta:h(\zeta)=0}\frac{\sqrt{2}}{\nu(\S)}\int_{h\in\S}\int_0^{\pi}\mu(h_t,\zeta_t)^2\,dt\;d\S.\right)=
\]
\[
\sqrt{2}\;{\rm E}_{g,h\in\S\times \S}\left(\int_{f\in L_{g,h}}\sum_{\zeta:f(\zeta)=0 }\mu(f,\zeta)^2 \right),
\]
where $L_{g,h}$ is the half--great circle in $\S$ containing $g,h$, starting at $g$ and going to $-g$ (we have to remove from this argument the case $h=-g$ but this is unimportant for integration purposes.) Note that we can define a measure and more generally a concept of integral in $\S$ as follows: given any measurable function $q:\S\rightarrow[0,\infty)$, its integral is
\begin{equation}\label{eq:mi}
{\rm E}_{g,h\in\S\times \S}\left(\int_{f\in L_{g,h}}q(f)\right).
\end{equation}
Now, this last formula describes an invariant (with respect to the group of symmetries of $\S$, that can be identified with the unitary group of size $N+1$ or with the orthogonal group of size $2N+2$) measure in $\S$ and is thus equal to a multiple of the usual measure in $\S$. In words, averaging over $\S$ or over great circles in $\S$ is the same, up to a constant. The constant is easy to compute by considering the constant function $q\equiv 1$. What we get is:
\[
{\rm E}_{g\in\S}\left( \sum_{\zeta:g(\zeta)=0}\mathcal{A}_1(g,\zeta)\right)\leq \frac{\pi}{\sqrt{2}}{\rm E}_{h\in\S}\left( \sum_{\zeta:h(\zeta)=0}\mu(h,\zeta)^2\right).
\]
After this argument is made rigorous, we have (see \cite{BePa:FoCM,BePa:JAMS} for earlier versions of the following result:)
\begin{theorem}[\cite{BePa:JFoCM}]\label{th:bp1}
Let $g\in\S$ be chosen at random with the uniform distribution, and let $\zeta$ be chosen at random, with the uniform (discrete) distribution among the roots of $g$. Then, the expected value of $\mathcal{A}_1(g,\zeta)$ is at most $\frac{\pi}{\sqrt{2}} nN$. In particular, for such a randomly chosen pair $(g,\zeta)$, with probability at least $1/2$ we have $\mathcal{A}_1(g,\zeta)\leq \sqrt{2}\pi nN$, that is, $(g,\zeta)$ is a good starting pair\footnote{Note that we are computing there the average of $\mathcal{A}_1$ not that of the integral of $\mu^2$ as in \cite{BePa:JFoCM}. From (\ref{eq:A1mu2}), the constant $\sqrt{2}$ has to be added to the formula in \cite{BePa:JFoCM} in this context.}.
\end{theorem}
The previous result would be useless for describing an algorithm (because choosing a random zero of a randomly chosen $g\in\S$ might be a difficult problem) without the following one.
\begin{theorem}[\cite{BePa:JFoCM}]\label{th:bp2}
The process of choosing a random $g\in\S$ and a random zero $\zeta$ of $g$ can be emulated by a simple linear algebra procedure.
\end{theorem}
The details of the linear algebra procedure of Theorem \ref{th:bp2} require the introduction of too much notation. We just describe the process in words: one has to choose a random $n\times(n+1)$ matrix $M$ with complex entries, compute its kernel (a projective point $\zeta\in\Pn$) and consider the system $g\in\S$ that has $\zeta$ as a zero and whose linear part is given by $M$. A random higher--degree term has to be added to $g$, and then linear and higher--degree terms must be correctly weighted. This whole process has running time polynomial in $N$. We thus have:
\begin{corollary}\label{cor:BP}
The linear homotopy algorithm with the starting pair obtained as in Theorem \ref{th:bp2} has average complexity\footnote{We use here the $\tilde{O}(X)$ notation: this is the same as $O(X\log(X)^c)$ for some constant $c$, that is logarithmic factors are cleaned up to make formulas look prettier.} $\tilde{O}(N^2)$.
\end{corollary}

The word ``average'' in Corollary \ref{cor:BP} must be understood as follows. For an input system $h$, let $T(h)$ be the average running time of the linear homotopy algorithm, when $(g,\zeta)$ is randomly chosen following the procedure of Theorem \ref{th:bp2}. Then, the expected value of $T(h)$ for random $h$ is $\tilde{O}(N^2)$. This kind of algorithm is called { Average Las Vegas}, the ``Las Vegas'' term coming from the fact that a random choice has to be done. The user of the algorithm plays the role of a Las Vegas casino, not of a Las Vegas gambler: the chances of winning (i.e. getting a fast answer to our problem) are much higher than those of loosing (i.e. waiting for a long time before getting an answer.)

Some of the higher moments of $\mathcal{A}_1(g,\zeta)$ have also been proved to be small. For example, the second moment (thus, also the variance) of $\mathcal{A}_1(g,\zeta)$ is polynomial in $N$, as the following result shows:
\begin{theorem}[\cite{BeltranShub2009}]\label{th:finitevariance}
Let $2\leq k<3$. Then, the expectation of $\mathcal{A}_1(g,\zeta)^k$ satisfies
\[
{\rm E}\left(\mathcal{A}_1(g,\zeta)^k\right)<\infty.
\]
Moreover, let $2\leq k<3-\frac{1}{2\ln\CD}$. Then,
the expectation ${\rm E}\left(\mathcal{A}_1(g,\zeta)^k\right)$ satisfies,
\[
{\rm E}\left(\mathcal{A}_1(g,\zeta)^k\right)\leq 2^{2k+k/2+4}\;e\;\pi^{k}n^{3k-4}N^2\CD^{4k-8}\ln\CD.
\]
In particular, ${\rm E}(\mathcal{A}_1(g,\zeta)^2)\leq 512 e \pi^2n^2N^2\ln\CD$.
\end{theorem}

We have been concentrating on finding one zero of a polynomial system. But we could find $k$ zeros $0\leq k\leq \CD$ by choosing $k$ different random initial pairs using Theorem \ref{th:bp2}. This process is known from \cite{BePa:JFoCM} to ouput every zero of the goal system $h$ with the same probability $1/\CD$, if $h\not\in\Sigma$. Another option is to choose some initial system $g$ which has $k$ known zeros, and simultaneously continuing the $k$ homotopy paths with the algorithm of Theorem \ref{th:BezVI}. In the case of finding all zeros the sum of the number of steps to follow each path, is by Theorem \ref{th:BezVI} and (\ref{eq:muboundszetadot}) bounded above by a constant times
\[
d^{3/2}\int_{0}^{d_R(g,h)}\sum_{\zeta_t:h_t(\zeta_t)=0}\mu(h_t,\zeta_t)^2\,dt.
\]
So for the great circle homotopies we have been discussing an analogous of Theorem \ref{th:bp1} holds:
\begin{theorem}[\cite{BePa:JFoCM}]\label{th:allzeros}
Let $g\in\S$ be chosen at random with the uniform distribution. Then, the expected value of $\int_{0}^{d_R(g,h)}\sum_{\zeta_t:h_t(\zeta_t)=0}\mu(h_t,\zeta_t)^2\,dt$ is at most $\frac{\pi}{2} nN\CD$. In particular, for such a randomly chosen $g$, with probability at least $1/2$ we have $\int_{0}^{d_R(g,h)}\sum_{\zeta_t:h_t(\zeta_t)=0}\mu(h_t,\zeta_t)^2\,dt\leq \pi nN\CD$, that is, the linear homotopy for finding all zeros starting at $g$ takes at most a constant times $d^{3/2}nN\CD$ steps to output all zeros of $h$, on the average.
\end{theorem}
Note that in general, one cannot write down all the $\CD$ zeros of $g$ to begin with, so Theorem \ref{th:allzeros} does not immediately yield a practical algorithm.

We point out that, even for the case $n=1$, no explicit descriptions of pairs $(g,\zeta)$ satisfying $\mathcal{A}_1(g,\zeta)\leq d^{O(1)}$ are known. Of course, no explicit polynomial $g\in\S$ is known in that case satisfying the claim of Theorem \ref{th:allzeros}. An attempt to determine such a polynomial lead to some progress in the understanding of elliptic Fekete points, see Section \ref{sec:fekete}.

\subsection{The roots of unity combined with a method of Renegar: a quasi--polynomial time deterministic algorithm for problem (\ref{prob:17})}
One can also ask for an algorithm for Problem (\ref{prob:17}) which does not rely on random choices (a { deterministic} algorithm). The search of a deterministic algorithm with polynomial running time for Problem (\ref{prob:17}) is still open, but a quasi--polynomial algorithm is known since \cite{BurgisserCucker2009}.

This algorithm is actually a combination of two: on one hand, we consider the initial pair
\begin{equation}\label{eq:BCpair}
g=
\begin{cases}
\frac{1}{\sqrt{2n}}(x_0^{d_1}-x_1^{d_1})\\
\vdots\\
\frac{1}{\sqrt{2n}}(x_0^{d_1}-x_n^{d_n})\\
\end{cases}
,\quad \zeta=(1,\ldots,1)
\end{equation}
Then, we have:
\begin{theorem}[\cite{BurgisserCucker2009}]\label{th:bc1}
The projective Newton--based homotopy method with initial pair (\ref{eq:BCpair}) has average running time polynomial in $N$ and $n^d$, where $d$ is the maximum of the degrees.
\end{theorem}

Theorem \ref{th:bc1} is a consequence of the following stronger result:
\begin{theorem}[\cite{BurgisserCucker2009}]\label{th:bc1b}
The projective Newton--based homotopy method with initial pair $(g,\zeta)\in\CV$ has average running time polynomial in $N$ and in $\max\{\mu(g,\eta):g(\eta)=0\}$.
\end{theorem}
Theorem \ref{th:bc1} follows from Theorem \ref{th:bc1b} and the fact that $\mu(g,\eta)\leq2(n+1)^d$ for $g$ given by (\ref{eq:BCpair}) for every zero $\eta$ of $g$.

For small (say, bounded) values of $d$, the quantity $n^d$ is polynomial in $n$ and thus polynomial in $N$, but for big values of $d$ the quantity $n^d$ is not bounded by a polynomial in $N$, and thus Theorem \ref{th:bc1} does not claim the existence of a polynomial running time algorithm. However, it turns out that there is a previously known algorithm, based on the factorization of the $u$--resultant, that has exponential running time for small degrees, but polynomial running time for high degrees (this may seem contradictory, but it is not: when the degrees are very high, the input size is big, and thus bounding the running time by a polynomial in the input size is sometimes possible in this case.) More precisely:

\begin{theorem}[\cite{Renegar89,BurgisserCucker2009}]\label{th:ren}
There is an algorithm with average running time polynomial in $N$ and $\CD$ that, on input $h\in\PHd\setminus\Sigma$, outputs an approximate zero associated to every single exact zero of $h$.
\end{theorem}

Note that $\CD$ is usually exponential in $n$, but as suggested above, if the degrees are very high compared to $n$, then $\CD$ can be bounded above by a polynomial in the input size $N$ and thus the algorithm of Theorem \ref{th:ren} becomes a polynomial running time algorithm.

An appropriate combination of theorems \ref{th:bc1} and  \ref{th:ren}, using the homotopy method of Theorem \ref{th:bc1} for moderately low degrees and the symbolic--numeric method of Theorem \ref{th:ren} for moderately high degrees (or, simply, running both algorithms for every input and stopping whenever one of them finishes,) turns out to be quasipolynomial for every choice of $n$ and $(d)$ . Indeed:
\begin{theorem}[\cite{BurgisserCucker2009}]\label{th:bc2}
The average (for random $h\in\S$) running time of the following procedure is $O(N^{\log\log N})$: on every input $h\in\PHd\setminus\Sigma$, run simultaneously the algorithms of theorems \ref{th:bc1} and  \ref{th:ren}, stopping the computation whenever one of the two algorithms gives an output.
\end{theorem}
Note that the running time of this algorithm is thus quasi--polynomial in $N$. Moreover, the algorithm is deterministic because it does not involve random choices.

\subsection{Homotopy paths based in the evaluation at one point}
Another approach to construct homotopies was considered in \cite{Smale81} and generalized in \cite{ArmentanoShub2012}. Given $h\in\CHd$ and $\zeta\in\Pn$, consider $g=h-\hat{h}_\zeta$, where $\hat{h}_\zeta\in\CHd$ is defined as
\[
\hat{h}_\zeta(z)=Diag\left(\frac{\langle z,\zeta\rangle^{d_i}}{\langle\zeta,\zeta\rangle^{d_i}}\right)h(\zeta).
\]
Then, $g(\zeta)=0$. So, we consider the homotopy $h_t=(1-t)g+th=h-(1-t)\hat{h}_\zeta$. We continue the zero $\zeta$ from $h_0=g$ to $h_1=h$. For any fixed $\zeta$, for example $\zeta=e_0=(1,0,\ldots,0)$, the homotopy may be continued for almost all $h\in\CHd$. Let
\[
K(h,\zeta)=\text{ number of steps sufficient to continue $\zeta$ to a zero of $h$}.
\]
Then,
\begin{theorem}[\cite{ArmentanoShub2012}]\label{th:armentanoshub}
\[
{\rm E}_{h\in\CHd}(K(h))\leq\frac{Cd^{3/2}\Gamma(n+1)2^{n-1}}{(2\pi)^N\pi^n}\int_{h\in\CHd}\left(\sum_{\eta:h(\eta)=0}\frac{\mu(h,\eta)^2}{\|h\|^2}\Theta(h,\eta)\right)e^{-\|h\|^2/2}\,dh,
\]
where
\[
\Theta(h,\eta)=\int_{\zeta\in B(h,\eta)}\frac{(\|h\|^2-T^2)1/2}{T^{2n-1}}\Gamma(T^2/2,n)e^{T^2/2}\,d\zeta,
\]
\[
T=\|Diag(\|\zeta\|^{-d_i}))h(\zeta)\|,
\]
and $\Gamma(\alpha,n)=\int_\alpha^{+\infty}t^{n-1}e^{-t}\,dt$ is the incomplete gamma function.
\end{theorem}
In Theorem \ref{th:armentanoshub}, $B(h,\eta)$ is the basin of $\eta$, which we now define. Suppose $\eta$ is a non--degenerate zero of $h\in\CHd$. We define the basin of $\eta$, $B(h,\eta)$, as those $\zeta\in\Pn$ such that the zero $\zeta$ of $g=h-\hat{h}_\zeta$ continues to $\eta$ for the homotopy $h_t=(1-t)g+th$. We observe that the basins are open sets.

Not much is known about ${\rm E}(K)$. Might it be polynomial in $N$? Is it even finite? See \cite{ArmentanoShub2012} for precise questions and motivations. Here is one:
\begin{problem}
For $h\in\CHd\setminus\Sigma$, are the volumes of all the basins of $h$ equal?
\end{problem}

\section{The condition Lipschitz--Riemannian structure}
Let us know turn our sight back to (\ref{eq:condmet}). If we drop the condition number $\mu(h_t,\zeta_t)$ from that formula, we get
\[
L=\int_a^b\|\dot h_t,\dot\zeta_t)\|\,dt,
\]
that is simply the length of the path $(h_t,\zeta_t)$ in the solution variety $\CV$, taking on $\CV$ the natural metric: the one inherited from that of the product $\PHd\times\Pn$. The formula in (\ref{eq:condmet}) can now be seen under a geometrical perspective: $L_\kappa$ is just the length of the path $(h_t,\zeta_t)$ when $\CV$ is endowed with the conformal metric obtained by multiplying the natural one by the square of the condition number. Note that this new metric is only defined on $\CW=\CV\setminus\Sigma'$. We call this new metric the { condition metric} in $\CW$. This justifies the name { condition length} we have given to $L_\kappa$. Theorem \ref{th:BezVI} now reads simply as follows: the complexity of following a homotopy path $(h_t,\zeta_t)$ is at most a small constant $cd^{3/2}$ times the length of $(h_t,\zeta_t)$ in the condition metric. This makes the condition metric an interesting object of study: which are the theoretical properties of that metric? given $p,q\in\CW$, what is the condition length of the shortest path joining $p$ and $q$?

The first thing to point out is that $\mu$ is not a $C^1$ function, as it involves a matrix operator norm. However, $\mu$ is locally Lipschitz. Thus, the condition metric is not a Riemannian metric (usually, one demands smoothness or at least $C^1$ for Riemannian metrics,) but rather we may call it a Lipschitz--Riemannian structure. This departs from the topic of most available books and papers dealing with geometry of manifolds, but there are still some things we can say. It is convenient to take a tour to a slightly more general kind of problems; that's the reason for the following section.

\subsection{Conformal Lipschitz--Riemann structures and self--convexity}
Let $\M$ be a finite--dimensional Riemannian manifold, that is a smooth manifold with a smoothly varying inner product defined at the tangent space to each point $x\in\M$, let us denote it $\langle \cdot,\cdot\rangle_x$. Let $\alpha:\M\rightarrow[0,\infty)$ be a Lipschitz function that is there exists some constant $K\geq0$ such that
\[
|\alpha(x)-\alpha(y)|\leq K d_R(x,y),\quad \forall\,x,y\in\M,
\]
where $d_R(x,y)$ is the Riemannian distance from $x$ to $y$. Then, consider on each point $x\in\M$ the inner product $\langle\cdot,\cdot\rangle_{\alpha,x}=\alpha(x)\langle\cdot,\cdot\rangle_x$. Note that this need no longer be smoothly varying with $x$, for $\alpha(x)$ is just Liptschitz. We call such a structure a (conformal) Lipschitz--Riemannian structure in $\CM$, and call it the $\alpha$--structure. The condition length of a $C^1$ path $\gamma(t)\subseteq\M$, $a\leq t\leq b$, is just
\[
L_\alpha(\gamma)=\int_a^b\|\gamma(t)\|_{\alpha,\gamma(t)}\,dt=\int_a^b\langle\gamma(t),\gamma(t)\rangle_{\alpha,\gamma(t)}^{1/2}\,dt
\]
The distance between any to points $p,q\in\M$ in this $\alpha$--structure is defined as
\begin{equation}\label{eq:dist}
d_\alpha(p,q)=\inf_{\gamma(t)\subseteq\CV}L_\alpha(\gamma),\quad p,q\in\M,
\end{equation}
where the infimum is over all $C^1$ paths with $\gamma(0)=p,\gamma(1)=q$.

A path $\gamma(t)$, $a\leq t\leq b$ is called a { minimizing geodesic} if $L_\alpha(\gamma)=d_\alpha(\gamma(a),\gamma(b))$ and $\|\dot\gamma(t)\|_{\alpha,\gamma(t)}\equiv1$, that is if it minimizes the length of curves joining its extremal points and if it is parametrized by arc--length. Then, a curve $\gamma(t)\subseteq\M$, for $t$ in some (possibly unbounded) interval $I$ is a { geodesic} if it is locally minimizing, namely if for every $t$ in the interior of $I$ there is some interval $(a,b)\subseteq I$ containing $t$ and such that $\gamma\mid_{[a,b]}$ is a minimizing geodesic.

Each connected component of the set $\M$ with the metric given by $d_\alpha$ is a path metric space, and it is locally compact because $\M$ is a smooth finite--dimensional manifold. We are in a position to use Gromov's version of the classical Hopf--Rinow theorem \cite[Th.1.10]{gro}, and we have:
\begin{theorem}\label{th:gromov}
Let $\M$ and $\alpha$ be as in the discussion above. Assume additionally that $\M$ is connected and that $(\M,d_\alpha)$ is a complete metric space. Then:
\begin{itemize}
\item each closed, bounded subset is compact,
\item each pair of points can be joined by a minimizing geodesic.
\end{itemize}
\end{theorem}
Theorem \ref{th:gromov} gives us sufficient conditions for conformal Lipschitz--Riemannian structures to be ``well defined'' in the sense that the infimum of (\ref{eq:dist}) becomes a minimum. We can go further:
\begin{theorem}[\cite{BeltranDedieuMalajovichShub2}]\label{th:BDMSa}
In the notation above, any geodesic is of class $C^{1+Lip}$, that is it is $C^1$ and has a Lipschitz derivative.
\end{theorem}
See \cite{BoitoDedieu} for an early version of Theorem \ref{th:BDMSa} and for experiments related to this problem.

One often thinks of the function $\alpha$ as some kind of ``squared inverse of the distance to a { bad} set'', so for each connected component of $\M$ the set $(\M,d_\alpha)$ will actually be complete.

A natural property to ask about is the following: given $p,q\in\M$, and given a geodesic $\gamma(t)$ such that $\gamma(a)=p$, $\gamma(b)=q$, does $\alpha$ attain its maximum on $\gamma$ in the extremes? Namely, if we think on $\alpha$ as some kind of squared inverse to a bad set, do we have to get closer to the bad set than what we are in the extremes?

\begin{example}\label{ex:pp}
A model to think of is Poincar\'e half--plane with the metric given by the usual scalar product in $\R^2\cap\{y>0\}$, multiplied by $y^{-2}$. Geodesics then become just portions of vertical lines or half--circles with center at the axis $y=0$. It is clear that, to join any two points, the geodesic does not need to become closer to the { bad} set $\{y=0\}$.
\end{example}

We can ask for more: we say that $\alpha$ is self--convex (an abbreviation for self--log--convex) if for any geodesic $\gamma(t)$, the following is a convex function:
\[
t\mapsto\log(\alpha(\gamma(t))).
\]
Note that this condition is stronger than just asking for $t\mapsto \alpha(\gamma(t))$ to be convex, and thus stronger than asking for the maximum of $\alpha$ on $\gamma$ to be at the extremal points.

\subsection{Convexity properties of the condition number}
We have the following result:
\begin{theorem}[\cite{BeltranDedieuMalajovichShub}]\label{th:BDMS1a}
Let $k\geq1$ and let $\mathrm{N}\subseteq\R^k$ be a $C^2$ submanifold without boundary of $\R^2$. Let $U\subseteq\R^n\setminus\mathrm{N}$ be the biggest open set all of whose points have a unique closest point in $\mathrm{N}$. Then, the function $\alpha(x)=distance(x,\mathrm{N})^{-2}$ is self--convex in $U$.
\end{theorem}
Note that Theorem \ref{th:BDMS1a} is a more general version of Example \ref{ex:pp}, where the horizontal line $\{y=0\}$ is changed to a submanifold $\mathrm{N}$.

A well--known result usually attributed to Eckart and Young \cite{EcYo36} and to Schmidt and Mirsky (see \cite{StSu90}) relates the usual condition number of a full rank rectangular matrix to the inverse distance to the set of rank--deficient matrices:
\begin{theorem}[Condition Number Theorem of linear algebra]\label{th:cnt}
Let $A\in\C^{mn}$ be a $m\times n$ matrix for some $1\leq m\leq n$. Let $\sigma_1(A),\ldots,\sigma_m(A)$ be its singular values. Then,
\[
\sigma_m(A)=distance(A,\{\text{rank--deficient matrices}\}).
\]
In particular, in the case of square maximal rank matrices, this we can rewrite this as $\|A^{-1}\|=distance(A,\{\text{rank--defficient matrices}\})^{-1}$, that is the (unscaled) condition number $\|A^{-1}\|$ equals the inverse of the distance from $A$ to the set of singular matrices. We more generally call $\sigma_m^{-1}(A)$ the unscaled condition number of a (possibly rectangular) full--rank matrix $A$.
\end{theorem}

One feels tempted to conclude from theorems \ref{th:BDMS1a} and \ref{th:cnt} that the function sending a full--rank complex matrix $A$ to the squared inverse of its smallest singular value (i.e. to the square of its unscaled condition number) should be self--convex. Indeed, one cannot apply Theorem \ref{th:BDMS1a} because the set of rank--deficient matrices is { not} a $C^2$ manifold, and because the distance to it is for many matrices (more precisely: whenever the multiplicity of the smallest singular value is greater than $1$) not attained in a single point. It takes a considerable effort to prove that the result is still true:
\begin{theorem}[\cite{BeltranDedieuMalajovichShub2}]\label{th:BDMS2b}
The function defined in the space of full--rank $m\times n$ matrices, $1\leq m\leq n$, as the squared inverse of the unscaled condition number, is self--convex.
\end{theorem}
Note that this implies that, given any two complex matrices $A,B$ of size $m\times n$, and given any geodesic $\gamma(t)$, $a\leq t\leq b$ in the $\alpha$--structure defined in
\[
\C^{mn}\setminus\{\text{ rank--deficient matrices}\}
\]
 by $\alpha(C)=\sigma_m(C)^{-2}$ such that $\gamma(a)=A$, $\gamma(b)=B$, the maximum of $\alpha$ along $\gamma$ is $\alpha(A)$ or $\alpha(B)$.

Note that, if a similar result could be stated for the $\alpha$--structure defined by $(h,\zeta)\mapsto\mu(h,\zeta)^2$ in $\CW$, we would have quite a nice description of how geodesics in the condition metric of $\CW$ are. Proving this is still an open problem:
\begin{problem}\label{prob:musc}
Prove or disprove $\mu^2$ is a self--convex function in $\CW$.
\end{problem}
Note that from Theorem \ref{th:condnumber}, the function $\mu^2$ is not exactly the squared inverse of the distance to a submanifold, but it is still something similar to that. This makes it plausible to believe that Problem \ref{prob:musc} has an affirmative answer. A partial answer is known:
\begin{theorem}[\cite{BeltranDedieuMalajovichShub2}]\label{th:BDMS2c}
The function $h\mapsto\mu^2(h,e_0)$ defined in the set $\{h\in\PHd:h(e_0)=0\}$ is self--convex. Here, $e_0=(1,0,\ldots,0)$.
\end{theorem}
\section{Condition geodesics and the geometry of $\CW$}
Although we do not have an answer Problem \ref{prob:musc}, we can actually state some bounds that give clues on the properties of the geodesics in the condition structure in $\CW$. More precisely:

\begin{theorem}[\cite{BeltranShub2007}]\label{cor:claro}
For every two pairs $(h_1,\zeta_1), (h_2,\zeta_2)\in \CW$, there exists a curve $\gamma_t\subseteq \CW$ joining $(h_1,\zeta_1)$ and $(h_2,\zeta_2)$, and such that
\[
L_\kappa(\gamma_t)\leq 2cnd^{3/2}+2\sqrt{n}\ln\left(\frac{\mu(h_1,\zeta_1)\mu(h_2,\zeta_2)}{n}\right),
\]
$c$ a universal constant.
\end{theorem}
In the light of Theorem \ref{th:BezVI}, this means that if one can find geodesics in the condition structure in $\CW$, one would be able to follow these paths in very few steps: just logarithmic in the condition number of the starting pair and the goal pair.
\begin{corollary}\label{cor:main1}
A sufficient number of projective Newton steps to follow some path in $\CW$ starting at the pair $(g,e_0)$ of (\ref{eq:shsmconj}) to find an approximate zero associated to a solution $\zeta$ of a given system $h\in\PHd$ is
\[
cd^{3/2}\left(nd^{3/2}+\sqrt{n}\ln\left(\frac{\mu(h,\zeta)}{\sqrt{n}}\right)\right),
\]
$c$ a universal constant.
\end{corollary}
Note that only the logarithm of the condition number appears in Corollary \ref{cor:main1}. Thus, if one could find an easy way to describe condition geodesics in $\CW$, the average complexity of approximating them using Theorem \ref{th:BezVI} would involve just the expectation of the average of $\ln(\mu)$, not that of $\mu^2$ as in Theorem \ref{th:averagemu2}. As a consequence, the average number of steps needed by such an algorithm would be $O(nd^3\ln N)$. See \cite[Cor. 3]{BeltranShub2009} for a more detailed statement of this fact.
At this point we ask a rather naive, vague question:
\begin{problem}\label{prob:la}
May homotopy methods be useful in solving linear systems of equations? Might using geodesics help as in Corollary \ref{cor:main1} and the comments above?
\end{problem}
Large sparse systems are frequently solved by iterative methods and the condition number plays a role in the error estimates. So Problem (\ref{prob:la}) has some plausibility.
\begin{remark}
There is an exponential gap between the average number of steps needed by linear homotopy $O(d^{3/2}nN)$ and those promised by the condition geodesic--based homotopy (which stays at a theoretical level by now, because one cannot easily describe those geodesics). This exponential gap occurs frequently in theoretical computer science. For example NP--complete problems are solvable in simply exponential time but polynomial with a witness. The estimates for homotopies with condition geodesics may likely serve as a lower bound for what can be achieved. Also, properties of geodesics as we learn them can inform the design of homotopy algorithms.
\end{remark}

There is more we can say about the geometry (and topology) of $\CW$, by studying the { Frobenius condition number} in $W$, which is defined as follows:
\[
\tmu(h,\zeta)=\|h\|\;\|Dh(\zeta)^\dagger\;
Diag(\|\zeta\|^{d_i-1}d_i^{1/2})\|_F,\;\;\;\forall (h,\zeta)\in W,
\]
where $\|\cdot\|_F$ is Frobenius norm (i.e. $Trace(L^* L)^{1/2}$ where $L^*$ is the adjoint of $L$) and $^\dagger$ is
Moore-Penrose pseudoinverse.
\begin{remark}The Moore-Penrose pseudoinverse $L^\dagger:\mathbb{F}\rightarrow \mathbb{E}$ of a linear operator $L:\mathbb{E}
\rightarrow \mathbb{F}$ of finite dimensional Hilbert spaces is defined as the composition
\begin{equation}\label{eq:MoorePe}
L^\dagger=i_{\mathbb{E}}\circ (L\mid_{Ker(L)^\perp})^{-1}\circ\pi_{Image(L)},
\end{equation}
where $\pi_{Image(L)}$ is the orthogonal projection on image $L$, $Ker(L)^\perp$ is the orthogonal complement of the nullspace of $L$, and $i_{\mathbb{E}}$ is the inclusion. If $A$ is a $m\times(n+1)$ matrix and $A=UDV^*$ is a singular value decomposition of $A$, $D=Diag(\sigma_1,\ldots,\sigma_k,0,\ldots,0)$ then we can write
\begin{equation}\label{eq:moorepenrosesvd}
A^\dagger=VD^\dagger U^*,\qquad D^\dagger=Diag(\sigma_1^{-1},\ldots,\sigma_k^{-1},0,\ldots,0).
\end{equation}
\end{remark}
In \cite{BeltranShub2009b} we prove that $\tmu$ is an equivariant Morse function defined in $\CW$ with a unique orbit of minima given by the orbit $\mathcal{B}$ of the pair of (\ref{eq:shsmconj}) under the action of the unitary group $(U,(h,\zeta))\mapsto (h\circ U^*,U\zeta)$.

The function $\mathcal{A}_1(g,\zeta)$ or even its upper bound (up to a $\sqrt{2}$ factor) estimate \footnote{see (\ref{eq:A1mu2}).}
\[
\mathcal{B}_1(g,\zeta)=\frac{1}{\nu(\S)}\int_{h\in\S}\int_0^\pi\mu(h_t,\zeta_t)^2\,dt\,d\S
\]
is an average of $\mu^2$ in great circles. This remark motivates the following \begin{problem}
Is $\mathcal{A}_1(g,\zeta)$ or $\mathcal{B}_1(g,\zeta)$ also an equivariant Morse function whose only critical point set is a unique orbit of minima?
\end{problem}
If so, due to symmetry considerations, it is the orbit through the conjectured good starting point (\ref{eq:shsmconj}). Here, one may want to replace the condition number $\mu$ in the definition of $\mathcal{B}_1$ with a smooth version such as the Frobenius condition number. A positive solution to this problem solves our main problem: the conjectured good initial pair (\ref{eq:shsmconj}) is not only good but even best.

Because the Frobenius condition number is an equivariant Morse function, the homotopy groups of $\CW$ are equal to those of $\mathcal{B}$, that can be studied with standard tools from algebraic topology. In the case that $n>1$, for example, we get:
\begin{align*}
\pi_0(\CW)&=\{0\}\\
\pi_1(\CW)&=\Z/a\Z\\
\pi_2(\CW)&=\Z\\
\pi_3(\CW)&=\pi_k(\mathcal{SU}_{n+1})\;(k\geq 3),
\end{align*}
where $\mathcal{SU}_{n+1}$ is the set of special unitary matrices of size $n+1$, $a=gcd(n,d_1+\cdots+d_n-1)$ and $\Z/a\Z$ is the finite cyclic group of $a$ elements.

In particular, we see that if all the $d_i's$ are equal then $a=1$ and $\CW$ is simply connected; in particular, any curve can be continously deformed into a minimizing geodesic. See \cite{BeltranShub2009b} for more results concerning the geometry of $\CW$. We can also prove a lower bound similar to the upper bound of Theorem \ref{cor:claro}:
\begin{theorem}\label{th:lbgeod}
let $\alpha:[a,b]\rightarrow\CW$ be a $C^1$ curve. Then, its condition length is at least
\[
\frac{1}{d^{3/2}\sqrt{n+1}}\left|\ln\left(\frac{\mu(\alpha(a))}{\mu(\alpha(b))}\right)-\ln\sqrt{n+1}\right|.
\]
\end{theorem}
\begin{remark}
We have written Theorem \ref{th:lbgeod} using the condition metric as defined in this paper. The original result \cite[Prop. 11]{BeltranShub2009b} was written for the so--called smooth condition length, obtained by changing $\mu$ to $\tmu$ in the definition of the condition length. This change produces the $\sqrt{n+1}$ factors in Theorem \ref{th:lbgeod}.
\end{remark}

In his article \cite{Smale1990}, Smale suggests that the input size of an instante of a numerical analysis problem should be augmented by $\log W(y)$ where $W(y)$ is a weight function ``... to be chosen with much thought...'' and he suggests that `` the weight is to resemble the reciprocal of the distance to the set of ill--posed problems.'' That is the case here. The condition numbers we have been using are comparable to the distance to the ill--posed problems and figure in the cost estimates. It would be good to develop a theory of computation which incorporates the distance to ill--posedness, or condition number and distance to ill-posedness in case they may not be comparable, (and precision in the case of round--off error) more systematically so that a weight function will not require additional thought. For the case of linear programming Renegar \cite{Renegar1995SIAM} accomplished this. It is our main motivating example as well as the work we have described on polynomial systems. The book \cite{BurgisserCuckerBook} is the current state of the art. The geometry of the condition metric will to our mind intervene in the analysis. If floating point arithmetic is the model of arithmetic used then ill-posedness will include points where the output is zero as well as points where the output is not Lipschitz.

\section{The univariate case and elliptic Fekete points}\label{sec:fekete}
Let us now center our attention in the univariate case, that, once homogenized, turns out into the case of degree $d$ homogeneous polynomials in two variables. One can then compute explicitly:
\[
\mu(h,\zeta)=d^{1/2}\|\|h\|\|(Dh(\zeta)\mid_{\zeta^\perp})^{-1}\|\zeta\|^{d-1}.
\]
If we are given a univariate polynomial $f(x)$ and a complex zero $z$ of $f$, we can also use the following more direct (and equivalent) formula for $\mu(h,\zeta)$ where $h$ is the homogeneous counterpart of $f$ and $\zeta=(1,z)$ :
\[
\mu(h,\zeta)=\frac{d^{1/2}(1+|z|^2)^{\frac{d-2}{2}}}{|f'(z)|}\|h\|.
\]
It was noted in \cite{ShSm93} that the condition number is related to the classical problem of finding elliptic Fekete points, which we recall now in its computational form (see \cite{BeltranFocmvolume} for a survey on the state of art of this problem.)

Given $d$ different points $x_1,\ldots,x_d\in\R^3$, let $X=(x_1,\ldots,x_d)$ and
\[
\EF(X)=\EF(x_1,\ldots,x_d)=-\sum_{i<j}\log\|x_i-x_j\|
\]
be its { logarithmic potential}. Sometimes $\EF(X)$ is denoted by $\EF_0(X)$, $\EF(0,X)$ or $V_N(X)$. Let $S(1/2)$ be the Riemann sphere in $\R^3$, that is the sphere of radius $1/2$ centered at $(0,0,1/2)$, and let
\[
m_d=\min_{x_1,\ldots,x_d\in S(1/2)}\EF(x_1,\ldots,x_d)
\]
be the minimum value of $\EF$. A minimising $d$--tuple $X=(x_1,\ldots,x_d)$ is called a set of elliptic Fekete points \footnote{Such a $d$--tuple can also be defined as a set of $d$ points in the sphere which maximize the product of their mutual distances.}.

The computational problem of finding elliptic Fekete points is another of the problems in Smale's list \footnote{Smale thinks on points in the unit sphere, but we may think on points in the Riemann sphere, as the two problems are equivalent by sending $(a,b,c)\in S(1/2)$ to $2(a,b,c)-(0,0,1)$.}.

{ Smale's 7th problem \cite{Smale00}:\;}{ Can one find $X=(x_1,\ldots,x_d)$ such that
\begin{equation}\label{eq:smale7}
\EF(X)-m_d\leq c\log d,\qquad c\text{ a universal constant}.
\end{equation}
}

The first clue that this problem is hard comes from the fact that the value of $m_d$ is not known, even to $O(d)$. A general technique (valid for Riemannian manifolds) given by Elkies shows that
\[
m_d\geq \frac{d^2}{4}-\frac{d\log d}{4}+O(d).
\]
Wagner \cite{Wagner1989} used the stereographic projection and Hadamard's inequality to get another lower bound. His method was refined by Rakhmanov, Saff and Zhou \cite{RakhmanovSaffZhou1994}, who also proved an upper bound for $m_d$ using partitions of the sphere. The lower bound was subsequently improved upon by Dubickas and Brauchart \cite{Dubickas1996}, \cite{Brauchart2008}. The following result summarizes the best known bounds:
\begin{theorem}\label{th:rsz}
Let $C_d$ be defined \footnote{The result in the original sources is written for the unit sphere, we translate it here to the Riemann sphere.} by
\[
m_d=\frac{d^2}{4}-\frac{d\log d}{4}+C_dd.
\]
Then,
\[
-0.4375\leq\liminf_{d\mapsto\infty}C_d\leq\limsup_{d\mapsto\infty}C_d\leq -0.3700708...
\]
\end{theorem}

The relation of this problem to the condition number relies in the fact that sets of elliptic Fekete points are naturally ``well separated'', and are thus good candidates to be the zeros of a ``well--conditioned'' polynomial, that is a polynomial all of whose zeros have a small condition number.
In \cite{ShSm93} Shub and Smale proved the following relation between the condition number and elliptic Fekete points.
\begin{theorem}[\cite{ShSm93c}]\label{th:ShSm93c}
Let $\zeta_1,\ldots,\zeta_d\in\mathbb{P}(\C^2)$ be a set of projective points, and consider them as points in the Riemann sphere $S(1/2)$ with the usual identification $\mathbb{P}(\C^2)\equiv S(1/2)$. Let $h$ be a degree $d$ homogeneous polynomial such that its zeros are $\zeta_1,\ldots,\zeta_d$. Then,
\[
\max\{\mu(h,\zeta_i):1\leq i\leq d\}\leq \sqrt{d(d+1)}\,e^{\EF(\zeta_1,\ldots,\zeta_d)-m_d}.
\]
In particular, is $x_1,\ldots,x_d$ are a set of elliptic Fekete points, then
\[
\max\{\mu(h,\zeta_i):1\leq i\leq d\}\leq\sqrt{d(d+1)}.
\]
\end{theorem}
\begin{remark}
Let $\mathfrak{Re}$ and $\mathfrak{Im}$ be, respectively, the real and complex part of a complex number. Here is alternative, equivalent definition for $h$ and the $\zeta_i$. Instead of considering projective points in $\mathbb{P}(\C^2)$ we may just consider a set of complex numbers $z_1,\ldots,z_d\in\C$. Then, for $1\leq i\leq d$, we can define $\zeta_i\in\S$ as
\begin{equation}\label{eq:xizi}
\zeta_i=\left(\frac{\mathfrak{Re}(z_i)}{1+|z_i|^2},\frac{\mathfrak{Im}(z_i)}{1+|z_i|^2},\frac{1}{1+|z_i|^2}\right)^T\in S(1/2),\qquad 1\leq i\leq d,
\end{equation}
$f$ as the polynomial whose zeros are $z_1,\ldots,z_d$, and $h$ as the homogeneous counterpart of $f$.
\end{remark}

There exists no explicit known way of describing a sequence of polynomials satisfying $\max\{\mu(h,\zeta):h(\zeta)=0\}\leq d^c$, for any fixed constant $c$ and $d\geq1$. Theorem \ref{th:ShSm93c} implies that, if a $d$--tuple satisfying (\ref{eq:smale7}) can be described for any $d$, then such a sequence of polynomials can also be generated. From Theorem \ref{th:bc1b}, such $h's$ are good starting points for the linear homotopy method, both for finding one root and for finding all roots. So, solving the elliptic Fekete points problem solves the starting point problem for $n=1$. The reciprocal question is: does solving the starting point problem for $n=1$ help with the Fekete point problem?
\begin{problem}
Suppose $n=1$ and $g\in\S$ minimizes $\sum_{\zeta:g(\zeta)=0}\mu(g,\zeta)^2$. Do $\zeta_1,\ldots,\zeta_d$ (the zeros of $g$, seen as points in $S(1/2)$) solve Smale's 7--th problem?
\end{problem}

We have seen in Theorem \ref{th:averagemu2} that the condition number of $(h,\zeta)$ where $h$ is chosen at random and $\zeta$ is uniformly chosen at random among the zeros of $h$, grows polynomially in $d$. Then, Theorem \ref{th:ShSm93c} suggests that spherical points associated with zeros of random polynomials might produce small values of $\EF$. We can actually put some numbers to this idea. First, one can easily compute the average value of $\EF$ when $x_1,\ldots,x_d$ are chosen at random in $S(1/2)$, uniformly and independently with respect to the probability distribution induced by Lebesgue measure in $S(1/2)$:
\[
{\rm E}_{X\in S(1/2)^d}\EF(X)=\frac{d^2}{4}-\frac{d}{4}.
\]
By comparing this with Theorem \ref{th:rsz}, we can see that random choices of points in the sphere already produce pretty low values of the minimal energy. One can prove that random polynomials actually produce points which behave better with respect to $\EF$:
\begin{theorem}[\cite{ArmentanoBeltranShub2009}]\label{th:abs}
Let $n=1$ and $h\in\S$ be chosen at random w.r.t. the uniform distribution in $\S$. Let $\zeta_1,\ldots,\zeta_d\in S(1/2)$ be the zeros of $h$. Then, the expected value of $\EF(\zeta_1,\ldots,\zeta_d)$ equals
\[
\frac{d^2}{4}-\frac{d\log d}{4}-\frac{d}{4}.
\]
\end{theorem}
By comparing this with Theorem \ref{th:rsz}, we conclude that spherical points coming from zeros of random polynomials agree with the minimal value of $\EF$, to order $O(d)$. This result fits into a more general (yet, less precise) kind of result related to random sections on Riemann surfaces, see \cite{Zhong2008,ZelditchZhong2010}.

\section{The algebraic eigenvalue problem}
The double fibration scheme proposed in (\ref{eq:pi1pi2}) has been -- at least partly -- successfully used in other contexts. For example, in \cite{Diego} a similar projection scheme
\begin{equation}\label{eq:pi1pi2eig}
\begin{matrix}
&&\CV_{eig}&&&\hskip -2cm =\{(A,\lambda,v)\in\mathbb{P}(\C^{n^2+1})\times\P(\C^n):Av=\lambda v\}\\
&\pi_1\swarrow&&\searrow\pi_2&&\\
&&&&&\\
\PHd&&&&\Pn&
\end{matrix}
\end{equation}
was used to study the complexity of a homotopy--based eigenvalue algorithm, obtaining the following:
\begin{theorem}\label{th:diego}
An homotopy algorithm can be designed that continues an eigenvalue--eigenvector pair $(\lambda_0,v_0)$ of a $n\times n$ matrix $A_0$ to one $(\lambda_1,v_1)$ of another matrix $A_1$, the number of steps bounded above by
\[
c\int_0^1\|(\dot A,\dot \lambda,\dot v)\|\mu_{eig}(A,\lambda,v)\;dt,
\]
$c$ a universal constant. Here, $\mu_{eig}$ is the condition number \footnote{A quantity similar in spirit to the condition number $\mu$ for the polynomial system solving problem.} for the algebraic eigenvalue problem , defined as
\begin{equation}\label{eq:condlengtheig}
\mu_{eig}(A,\lambda,v)=\max\left\{1,\|A\|_F\|\pi_{v^\perp}(\lambda I_n-A)\mid_{v^\perp}^{-1}\|\right\},
\end{equation}
where $\|A\|_F=trace(A^*A)^{1/2}$ is the Frobenius norm of $A$.
\end{theorem}
Of course, we do not intend to summarize here the enormous amount of methods and papers dealing with the eigenvalue problem (see \cite{StSu90} for example). We just point out that there exists no proven polynomial--time algorithm for approximating eigenvalues (although different numerical methods achieve spectacular results in practice.) See \cite{PfrangDeiftMenon} for some statistics about the QR (and Toda) algorithms for symmetric matrices. We don't know a good reference for the more difficult general case. Unshifted QR is not the fast algorithm of choice. The QR algorithm with Francis double shift executed on upper Hermitian matrices should be the gold standard.
\begin{problem}
Does the QR algorithm with Francis double shift fail to attain convergence on an open subset of upper Hessenberg matrices?
\end{problem}
See \cite{BattersonSmilie1989} for open sets where Rayleigh quotient iteration fails, and \cite{Batterson1994} for a proof of convergence for normal matrices as well as a good introduction to the dynamics involved.

 Theorem \ref{th:diego} can probably be used in an analysis similar to that of Section \ref{sec:start} to complete a complexity analysis. Note that the integral in Theorem \ref{th:diego} is very similar in spirit to that of (\ref{eq:condmet}). This allows to introduce a condition metric in $\CV_{eig}$. Some of the results in previous sections can be adapted to this new case. For example, an analogous of Theorem \ref{cor:claro} holds (i.e. short geodesics exist,) see \cite{diegotesis}.

The eigenvalue problem and the problem of finding roots of a polynomial in one variable are, of course, connected. Given an $n\times n$ matrix $A$ we may compute the characteristic polynomial of $A$, $p(z)=det(zI-A)$ and then solve $p(z)$. The zeros of $p(z)$ are the eigenvalues of $A$. Trefethen and Bau \cite{TrefethenBau1997} write ``This algorithm is not only backward unstable but unstable and should not be used''. Indeed when presented with a univariate polynomial $p(z)$ to solve numerical linear algebra packages may convert the problem to an eigenvalue problem by considering the companion matrix of $p(z)$ and then solve the eigenvalue problem. If $p(z)=z^d+a_{d-1}z^{d-1}+\cdots+a_0$ the compainon matrix is
\[
\begin{pmatrix}
0&0&0&\cdots&0&-a_0\\
1&0&0&\cdots&0&-a_1\\
0&1&0&\cdots&0&-a_{2}\\
\vdots& &\ddots&\ddots&\vdots&\vdots\\
\vdots& &\ddots&\ddots&0&-a_{d-2}\\
0&\cdots&\cdots&0&1&-a_{d-1}
\end{pmatrix},
\]
which is already in upper Hessenberg form. So conceivably Francis double shifted QR may fail to converge on an open set of companion matrices?

Let us recall that the condition number of a polynomial and root is a property of the output map as a function of the input. So doesn't depend on the algorithms to solve the problem. This motivates the following
\begin{problem}
What might explain the experience of numerical analysts, relating the polynomial solving methods versus that of eigenvalue solving? Might the condition number of the eigenvalue problem have small average over the set of $n\times n$ matrices with a given characteristic polynomial?
\end{problem}
Finally, if we consider the equation $Av=\lambda v$ an equation in unknowns $\lambda$ and $v$, it is a set quadratic equation. So, we have $n$ quadratic equations in $n$ unknowns. By Bezout's theorem, after we homogenize, we expect $2^n$ roots counted with multiplicity. But there are only $n$ eigenvalues. In \cite{Diego,diegotesis} it is shown that the use of multihomogeneous Bez\'out theorem yields the correct zero count for this problem. Thus, a reasonable thing to do is to introduce a new variable $\alpha$ and consider the bilinear equation $A\alpha v=\lambda v$ which is bilinear in $(\alpha,\lambda)$ and $v$.
\begin{problem}[see \cite{DeSh00}] Prove an analogue of Theorem \ref{th:diego} in the general multihomogeneous setting.
\end{problem}

\appendix
\section{A model of computation for machines with round-off and input errors}
This section has been developed in discussions with Jean Pierre Dedieu and his colleagues Paola Boito and Guillaume Ch\`eze. We thank Felipe Cucker for helpful comments.

\subsection{Introduction}
During the second half of the 20th century, with the emergence of computers, algorithms have taken a spectacular place in mathematics, especially numerical algorithms (linear algebra, ode's, pde's, optimization), but also symbolic computation. In this context, complexity studies give a better understanding of the intrinsic difficulty of a problem, and describe the performance of algorithms which solve such problems. One can associate the classical Turing model to symbolic computation based on integer arithmetic, and the BSS model to scientific computation on real numbers. However this ideal picture suffers from an important defect. Scientific computation does not use the exact arithmetic of real numbers but floating-point numbers and a finite precision arithmetic. Thus, a numerical algorithm designed on real numbers and the same algorithm running in finite precision arithmetic give a priori two different results. Any numerical analysis undergraduate book has at least one chapter dealing with the precision of numerical computations. See for example \cite{TrefethenBau1997} or \cite{Higham}. Yet, there is no solid approach to the definition and study of a model of computation including this aspect, as well as the role that conditioning of problems should play in the complexity estimates.

 Besides linear algebra problems and iterative processes, a key point to bear in mind is that we sometimes use floating point computers to answer decision (i.e Yes/No) problems, as { is this matrix singular?} or { does this polynomial have a real zero?}. The first attempts to use round-off machines to study decision problems are \cite{CuckerDedieu}, and \cite{CuckerSmale}. The authors consider questions like: under which conditions is the decision taken by a BSS machine the same as the decision taken by the corresponding round-off machine? Or, under which conditions is the decision taken by a round-off machine for a given input the same as the decision taken by the BSS machine on a nearby input?

In these pages we point towards the development of a theory of finite precision computation via a description of round--off machines, size of an input, cost of a computation, single (resp. multiple) precision computations (a computation is ``single precision'' when a sufficient round-off unit $\delta$ to reach relative precision $u$ for any input $x$ in the considered range is proportional to $u$), finite precision computability and finite precision decidability. These concepts have to be related to the intrinsic characteristics of the problem: its condition number (the local Lipschitz constant of the solution map), and its posedness (the distance to ill-posed problems).

The model we propose is inspired by the BSS model but it stays close to real-life numerical computation. We prefer relative errors than absolute ones (this is the basis of the usual floating-point arithmetic.) We mention two papers of interest about the foundations of scientific computing, \cite{BravermanCook,Braverman}, with
a point of view different than ours.

\subsection{Round--off machine}
A round-off machine is an implementation of a BSS-machine accounting for input error and round-off error of computations. These errors may mimic a particular floating point arithmetic but are designed to be more general. In particular, they are not  tied down to a particular floating point model.  Let $\R^\infty$ be the disjoint union of the sets $\R^n$, $n\geq0$. For given $x\in\R^\infty$ we define $\|x\|=\max_i|x_i|$. A subset $U\subseteq\R^{\infty}$ is open if it is the disjoint union of $U_n$ with $U_n\subseteq\R^n$ an open set. For this topology, a mapping $f:\R^\infty\rightarrow\R$ is continuous iff each restriction $f_n=f\mid_{\R^n}$ is continuous.

A BSS machine $M$ is a directed graph with with several kinds of nodes including an input node, with input $x\in\R^{\infty}$, output nodes, computation nodes where rational functions are generally computed but here we restrict ourselves without loss of generality to the standard arithmetic operations, branching nodes (we branch on an inequality of the form $y\geq0$.) A machine is a decision machine when the output is $-1$ or $1$. The  halting set $\mathcal{H}$ of $M$ is the set of inputs giving rise to an output. We denote by $\mathcal{O} : \mathcal{H} \rightarrow\R^{\infty}$ the
output map. There are a few technical concepts (mainly the { input map} $I(x)$ and the { computing endomorphism} $H_M$) associated to $M$, the nonfamiliar reader may find formal definitions in \cite[Chapter 2]{BlCuShSm98}.

Given a BSS machine $M$ defined over the real numbers as above and $ 0 \leq \delta \leq 1 $, a round-off machine associated to $M$ and $\delta$ is another machine denoted $(M,\delta)$ which we will call a BSS round-off machine. The nodes and state space of $(M,\delta)$ are the same as for $M$. The input map $I_\delta$ of $(M,\delta)$ satisfies $ |I_\delta (x)_j - I(x)_j| < \delta |I(x)_j|$ that is to say the relative error of the input is less than $\delta$ for every coordinate $j$. The next node next state of $(M,\delta)$ at a computation node is similarly  a relative approximation of the next node next state map of $M$ and does not necessarily represent an actual computation. The $jth$ components of the next states satisfy $|H_{(M,\delta),state}(x)_j - H_{M,state}(x)_j| < \delta |H_{M,state}(x)_j|$.

Note that given $M$ and $\delta$, there are many machines satisfying this definition. For example, $M$ itself satisfies this definition for every $\delta$. The power of the definition is that certain claims will hold for { every} such a round--off machine, allowing us to use just the defining properties and not the particular structure of a given round--off machine.

\subsection{Computability}

In the sequel, we will only consider functions $f : \Omega \subseteq\R^\infty\rightarrow\R^\infty$ such that, for each $n$, the restriction $f_n$ of $f$ to $\Omega_n=\Omega\cap\R^n$ takes its values in $\R^m$ for an $m$ depending only on $n$.

Such a function is round-off computable when there exists  a BSS machine $M$  such that for any $x \in \Omega$ and any $0<\epsilon<1$, there exists a $\delta (x,\epsilon)$ such that any round--off machine $(M,\delta(x,\epsilon))$ outputs $\tilde{O}(x)$ with
\[
|\tilde{O}(x)_j-f(x)_j|\leq\epsilon |f(x)_j|,
\]
that is the output of $(M,\delta(x,\epsilon))$ is coordinatewise equal to $f(x)$ up to relative precision $\epsilon$.

\begin{example}
The function $f:\R^2\ra\R$, $f(x,y)=xy$ (we can let it be zero in $\R^\infty\setminus\R^2$) is round--off computable. Indeed, let $x,y\neq0$ and $0<\epsilon<1$. The output of a round--off machine $(M,\delta)$ associated to the natural BSS machine for computing $f(x,y)$ is a number
\[
z=xy(1+\delta_1)(1+\delta_2)(1+\delta_3),
\]
for some $\delta_1,\delta_2,\delta_3$ bounded in absolute value by $\delta$. It is useful to note the elementary inequality
\begin{align}\label{eq:elemnineq}
\left|\left(1+\frac{u}{n}\right)^n-1\right|\leq2u,\quad \forall\;0\leq|u|\leq1.
\end{align}
From this, we obviously have $|z-xy|\leq \epsilon|xy|$ by taking
\begin{align}\label{eq:x*y}
\delta((x,y),\epsilon)=\frac{\epsilon}{6},
\end{align}
The output of any round--off machine if $x=0$ or $y=0$ is clearly $0$, and hence the same value for $\epsilon$ of (\ref{eq:x*y}) suffices to satisfy the definition of computability.
\end{example}

\begin{example}\label{ex:prodxi}
The same argument proves that the function $f:\R^\infty\ra\R$ given by $f(x_1,\ldots,x_n)=x_1\cdots x_n$ is round--off computable (say, we compute first $x_1x_2$ then $x_1x_2x_3$ and so on) with
\begin{align}\label{eq:prodxi}
\delta((x,y),\epsilon)=\frac{\epsilon}{4n-2},
\end{align}
\end{example}

\begin{example}\label{ex:x+y}
A longer computation shows that the function $f:\{(x,y)\in\R^2:x+y\neq0\}\ra\R$, $f(x,y)=x+y$ (again, we let it be zero in $\R^\infty\setminus\R^2$) is also round--off computable. It suffices to take
\[
\delta((x,y),\epsilon)=\frac{\epsilon}{2\max\left(1,\left|\frac{x}{x+y}\right|,\left|\frac{y}{x+y}\right|\right)}.
\]
A more simple and still valid formula is
\begin{align}\label{eq:x+y}
\delta((x,y),\epsilon)=\frac{|x+y|}{3\sqrt{2}\sqrt{x^2+y^2}}\epsilon.
\end{align}
\end{example}

\begin{example}\label{ex:sumpos}
Let us now see that $f(x)=x_1+\ldots+x_n$ is round--off computable in the set $\Omega=\{x\in\R^\infty:x_i\geq0\;\forall i\}$. Indeed, let $0<\epsilon<1$ and let us consider the most simple BSS machine which computes first $x_1+x_2$, then adds $x_3$ and so on\footnote{This is not the algorithm of choice in practical programming but is sufficient for our purposes here.} A round--off machine with precision $\delta$ will produce, on input $x=(x_1,\ldots,x_n)$, a number
\[
x_1\left(\prod_{k=1}^{n}(1+\delta_1^{(k)})\right)+ x_2\left(\prod_{k=1}^{n}(1+\delta_2^{(k)})\right)+
x_3\left(\prod_{k=2}^{n}(1+\delta_3^{(k)})\right)+\cdots
+x_n\left(\prod_{k=n-1}^{n}(1+\delta_n^{(k)})\right),
\]
for some $\delta_i^{(k)}$ bounded in absolute value by $\delta$. This follows from the fact that, in addition to the input error on each coordinate, $x_1$ and $x_2$ go through $n-1$ additions (which generate $n+1$ errors), $x_3$ goes though $n-2$ additions and so on. Note that
\[
x_1(1-\delta)^n\leq x_1\left(\prod_{k=1}^{n}(1+\delta_1^{(k)})\right)\leq x_1(1+\delta)^n.
\]
Choosing $\delta=\alpha \epsilon/(2n)$, $0<\alpha\leq1$ and using (\ref{eq:elemnineq}) we conclude that
\[
\left|x_1\left(\prod_{k=1}^{n}(1+\delta_1^{(k)})\right)-x_1\right|\leq \alpha\epsilon x_1,
\]
and the same formula holds for $x_2,\ldots,x_n$. The output of a round--off machine thus satisfies
\[
\Tilde{O}(x)=\sum_{i=1}^n x_i(1+\alpha\epsilon_i),\quad 0\leq|\epsilon_i|\leq\epsilon.
\]
That is,
\[
\left|\Tilde{O}(x)-\sum_{i=1}^nx_i\right|=\sum_{i=1}^nx_i\alpha\epsilon_i\leq \sum_{i=1}^nx_i\alpha|\epsilon_i|\leq\alpha \epsilon\sum_{i=1}^nx_i,
\]
proving that $f(x)$ is round--off computable in that set (just take $\alpha=1$).

\end{example}
\begin{example}\label{ex:sumxi}
Let us now see that $f(x)=x_1+\ldots+x_n$ is round--off computable in the set $\Omega=\{x\in\R^\infty:\sum x_i\neq0\}$. We consider the BSS machine that first adds all the nonnegative numbers, call $a$ the result, then adds all the negative numbers, call $b$ the result, and then computes $a-b$. Let $0<\epsilon<1$. We note that from Example \ref{ex:sumpos} by choosing $\delta=\alpha\epsilon/(2n)$ (some $0<\alpha\leq1$) the round--off computation of the sum of positive (resp. negative) terms will be
\[
\tilde{a}=a(1+\alpha\epsilon_1),\qquad \tilde{b}=b(1+\alpha\epsilon_2),\text{ for some }0\leq|\epsilon_1|,|\epsilon_2|\leq\epsilon.
\]
From Example \ref{ex:x+y}, if we let
\[
\alpha=\frac{|a+b|}{3\sqrt{2}\sqrt{a^2+b^2}},
\]
that is if we let
\[
\delta(x,\epsilon)\leq\frac{|a+b|}{3\sqrt{2}\sqrt{a^2+b^2}}\frac{\epsilon}{2n},
\]
then $\tilde{O}(x)=\sum_ix_i$ up to relative precision $\epsilon$. Using that $a^2+b^2\leq n\sum x_i^2$ We can also use the formula
\begin{align}\label{eq:sumxi}
\delta(x,\epsilon)=\frac{\left|\sum_{i=1}^n x_i\right|}{6\sqrt{2}n^{3/2}\sqrt{\sum_{i=1}^n x_i^2}}\,\epsilon.
\end{align}

\end{example}
\begin{example}
Combining examples \ref{ex:prodxi} and \ref{ex:sumxi} we see that the evaluation map of any multivariate polynomial $p(x_1,\ldots,x_n)$ is round--off computable in the complement of its zero set (just compute first the monomials and them add all the results).
\end{example}
\subsection{Ill--conditioned instances, condition number, posedness}

Let us think of a function $f:\Omega\subseteq\R^\infty\ra\R^\infty$ as the solution map associated with some problem to be solved. The condition number associated with $f$ and $x$ measures the first--order (relative) componentwise or normwise variations of $f(x)$ in terms of the first--order (relative) variations of $x$. 

First assume that $f:\Omega \ra\R$, that is the function is real--valued. We say that $x\in\bar{\Omega}$ (the topological closure of $\Omega$) is well--conditioned when:
\begin{itemize}
\item Either $\|x\|\neq0$, and $f$ can be extended to a Lipschitz function defined in a neighborhood of $x$ in $\bar{\Omega}$ with $|f(x)|\neq0$. In that case we define the componentwise condition number by
\[
\kappa_f(x)=\limsup_{x'\mapsto x,x'\in \bar{\Omega}}\frac{\frac{|f(x')-f(x)|}{|f(x)|}}{\frac{\|x'-x\|}{\|x\|}},
\]
\item or $f$ is constant in a neighborhood of $x$ with $|f(x)|=0$. In this later case we define the condition number by $\kappa_f(x)=0$.
\end{itemize}
Otherwise, we say that $x\in\bar{\Omega}$ is ill--conditioned. The set of ill--conditioned instances is denoted by $\Sigma_f$, and for $x\in\Sigma_f$, we let $\kappa_f(x)=\infty$.

For a general $f:\Omega\ra\R^\infty$, we define
\[
\kappa_f(x)=\sup_j\kappa_{f_j}(x)\quad {(componentise condition number)}
\]
that is the condition number of $f$ is the supremum of the condition numbers of its coordinates. Sometimes it is more useful to consider the normwise condition number, that we denote by the same letter as the context should make clear which one is used on each problem:
\[
\kappa_f(x)=\limsup_{x'\mapsto x,x'\in \bar{\Omega}}\frac{\frac{\|f(x')-f(x)\|}{\|f(x)\|}}{\frac{\|x'-x\|}{\|x\|}} \quad {(normwise condition number)},
\]

We define the posedness of a problem instance $x$ with $\|x\|\neq0$ as the distance to ill--posed problems:
\[
\pi_f(x)=\frac{d(x,\Sigma_f)}{\|x\|}.
\]
The relation between condition number and posedness is an important but unclear problem. We may expect a relation of the type
\[
\pi_f(x)\approx\kappa_f(x)^{-1}
\]
(condition number theorem) or at least inequalities like
\[
C_1\pi_f(x)^{\rho_1}\leq \kappa_f(x)^{-1}\leq C_w\pi_f(x)^{\rho_2}
\]
for suitable positive constants $C_i,\rho_i$ (cf. Lojasiewicz's inequality.) To get such a relation ill--posed problems should correspond to infinite condition numbers, but this is not always the case. Consider for example the decision problem: Is $x^2 + y^2 \leq \Pi ?$ The problem is well conditioned except on the circle $x^2 + y^2 = \Pi$, but the distance to this circle determines the precision we need in the computation.

Let $K_f(x)=max(\kappa_f(x), \pi_f(x)^{-1})$.

\begin{example}
For $f(x)=x_1\cdots x_n$ defined in $\R^\infty$, it is easy to see that
\[
\kappa_f(x)=\sqrt{x_1^2+\cdots +x_n^2}\sqrt{\frac{1}{x_1^2}+\cdots+\frac{1}{x_n^2}},
\]
whenever $x_1,\ldots,x_n\neq0$. If $x_i=0$ for any $i$ then $\kappa_f(x)=\infty$.

On the other hand,
\[
\pi_f(x)=\frac{\min(|x_1|,\ldots,|x_n|)}{\sqrt{x_1^2+\cdots +x_n^2}}.
\]
Thus, we have
\[
\kappa_f(x)\leq \sqrt{x_1^2+\cdots +x_n^2}\sqrt{\frac{n}{\min(|x_1|,\ldots,|x_n|)^2}}=\sqrt{n}\pi_f(x)^{-1},
\]
and
\[
\kappa_f(x)\geq \sqrt{x_1^2+\cdots +x_n^2}\sqrt{\frac{1}{\min(|x_1|,\ldots,|x_n|)^2}}=\pi_f(x)^{-1}.
\]
Namely,
\[
\pi_f(x)^{-1}\leq\kappa_f(x)\leq\sqrt{n}\pi_f(x)^{-1}.
\]

\end{example}

\begin{example}
For $f(x)=x_1+\cdots+x_n$ defined in $\Omega=\{x\in\R^\infty:\sum x_i\neq0\}$, we have:
\begin{itemize}
\item For $(x,y)\in\Omega$, a simple computation shows that
\[
\kappa_f(x,y)=\frac{\sqrt{n}\sqrt{\sum x_i^2}}{|\sum x_i|}.
\]
\item For $x\in\partial\Omega$, that is $\sum x_i=0$, we have $\kappa_f(x)=\infty$.
\end{itemize}
Thus, we have
\[
\pi_f(x)=\frac{d(x,\{x:\sum x_i=0\})}{\sqrt{\sum x_i^2}}=\frac{|\sum x_i|}{\sqrt{n}\sqrt{\sum x_i^2}}=\kappa_f(x,y)^{-1}.
\]
Namely,
\begin{align}\label{eq:x+yK}
K_f(x,y)=\frac{\sqrt{n}\sqrt{\sum x_i^2}}{|\sum x_i|}.
\end{align}
\end{example}

\subsection{Single, multiple precision}

Let $f$ be a round--off computable function, and let $M$ be a BSS machine satisfying the definition of round--off computability above. This computation is single precision when for every $0<\epsilon<1$ there is a $\delta=\delta(\epsilon)$ such that any round--off machine $(M,\delta)$ attains relative precision $\epsilon$ for any input $x\in\Omega$, and such that
\begin{equation}\label{eq:sp}
\delta\geq \frac{c_0\epsilon}{K_f(x)^{c_2}\dim(x)^{c_3}}
\end{equation}
for some positive constants $c_0,c_2,c_3$. This computation is multiple precision when a $\delta$ exists such that
\begin{equation}\label{eq:mp}
\delta\geq \frac{c_0\epsilon^{c_1}}{K_f(x)^{c_2}\dim(x)^{c_3}},
\end{equation}
for some $c_1>1$. We say that the computation is strictly multiple precision when it is multiple precision but not single precision.

\begin{example}
The inductive, naive algorithm for computing the round--off computable function $f(x_1,\ldots,x_n)=x_1\cdots x_n$ defined in $\R^\infty$ is single precision, from (\ref{eq:prodxi}). The inductive, naive algorithm for computing the round--off computable function $f(x)=x_1+\cdots+x_n$ defined in $\{x\in\R^\infty:\sum x_i\neq0\}$ is single precision from (\ref{eq:sumxi}) and (\ref{eq:x+yK}).
\end{example}
\subsection{Size of an input}
In many practical problems, we want to specify an output precision $\epsilon$. From our definition of round--off computable function, given $x\in\Omega$ and $0<\epsilon<1$ a $\delta(x,\epsilon)$ will exist guaranteeing the desired precision, although it may be very hard to compute this $\delta$ in some cases. Moreover, from (\ref{eq:mp}), the number $K_f(x)$  will in general play a role in the value of $\delta(x,\epsilon)$ needed for any machine solving the problem. This dependence suggests that maybe the input should be considered as $(x,\epsilon)$ and not just as $x$. These thoughts justify our definition of the size of an input, which includes a term related to $\epsilon$ and another related to $K_f(x)$:
\begin{equation}\label{eq:inputsize}
ht(x)+|\log\epsilon|+\log(K_f(x)+1).
\end{equation}
where by $ht(x)$ we mean $1$ if $x=0$ and $1+ln(1+|ln(|x|)|)$ if $x\neq 0$. If $x \in \R^n$ then by $ht(x)$ we mean $n \max ht(x_j)$ for $j=1 \dots n$ and $x_j$ the components of $x$.

\subsection{Cost of a computation}
The cost of a computation on a round-off machine $(M,\delta)$ which outputs $\tilde{y}$ on input $x$ is
\[
T(x,\delta)\cdot\left(\max_iht(y^{(i)})+|\log\delta| \right),
\]
where $T(x,\delta)$ is the time for the computation to halt and
\[
x=y^{(0)},\ldots,y^{T(x,\delta)}=\tilde{y}
\]
are the different vectors computed by $(M,\delta)$ on input $x$.

We say that a function $f:\Omega\ra\R^\infty$ is polynomial cost computable if for every  $x\in\Omega$ and $0<\epsilon<1$ there is exists $\delta(x,\epsilon)$  such that $(M,\delta(x,\epsilon))$ computes $\tilde{y}$ which equals $f(x)$  to relative error $\epsilon$ with cost polynomially bounded by the input size (\ref{eq:inputsize}).

The most important cases of polynomial cost computability will be in the cases where we restrict the space of functions to single (multiple) precision functions, for example in the case of single precision to the definition of polynomial cost we add the restriction that $\delta(x,\epsilon)$ must satisfy (\ref{eq:sp}). These two possibilities (single or multiple precision) will give us two theories, both of which deserve to be worked out.

Now that we have the notion of polynomial cost the classes $P$ and $NP$ may be defined and the problem: Does $P=NP?$ stated.


\providecommand{\bysame}{\leavevmode\hbox
to3em{\hrulefill}\thinspace}
\providecommand{\MR}{\relax\ifhmode\unskip\space\fi MR }
\providecommand{\MRhref}[2]{%
  \href{http://www.ams.org/mathscinet-getitem?mr=#1}{#2}
} \providecommand{\href}[2]{#2}

\end{document}